\documentclass{amsart}
\usepackage{amssymb}
\usepackage{amscd}
\usepackage{mathrsfs}
\usepackage{verbatim}

\theoremstyle{plain}
\newtheorem{theorem}{Theorem}[section]
\theoremstyle{remark}
\newtheorem{remark}[theorem]{Remark}
\newtheorem{example}[theorem]{Example}
\theoremstyle{plain}
\newtheorem{corollary}[theorem]{Corollary}
\newtheorem{lemma}[theorem]{Lemma}
\newtheorem{proposition}[theorem]{Proposition}
\newtheorem{definition}[theorem]{Definition}

\numberwithin{equation}{section}

\font\Bbbten=msbm10 
\font\Bbbseven=msbm7
\font\Bbbfive=msbm5 
\newfam\Bbbfam

  \textfont\Bbbfam=\Bbbten
  \scriptfont\Bbbfam=\Bbbseven
  \scriptscriptfont\Bbbfam=\Bbbfive
\def\R{{\mathbb R}\,}

\def\C{{\mathbb C}\,}
\def\E{{\mathbb E}\,}

\def \qed{\hfill\hbox{\hskip 6pt\vrule width6pt height7pt
depth1pt  \hskip1pt}\bigskip}

\renewcommand{\O}{\Omega}
\newcommand{\V}{\overline V}
\newcommand{\limn}{\lim_{n\to\infty}}
\renewcommand{\i}{i_\infty}
\newcommand{\nn}{|\!|\!|}
\newcommand{\e}{\varepsilon}
\newcommand{\z}{^{\odot{}}}
\renewcommand{\o}{\omega}
\renewcommand{\l}{\lambda}
\renewcommand{\Im}{{\rm Im}}
\renewcommand{\Re}{{\rm Re}}
\newcommand{\n}{\Vert}
\newcommand{\embed}{\hookrightarrow}
\newcommand{\s}{^*{}}
\newcommand{\D}{{\mathscr D}}
\newcommand{\lb}{\langle}
\newcommand{\rb}{\rangle}
\renewcommand{\H}{{H_\infty}}

\renewcommand{\SS}{{\bf S}}
\newcommand{\WW}{{\bf W}}
\newcommand{\TT}{{\bf T}}
\newcommand{\XX}{{\bf X}}
\newcommand{\PP}{{\bf P}}
\newcommand{\ZZ}{{\bf Z}}
\renewcommand{\a}{\alpha}

\newcommand{\Q}{Q_\infty}
\renewcommand{\L}{\mathscr L}

\newcommand{\HQinfty}{{\bf H}Q_\infty}
\newcommand{\Hmu}{{\bf H}\mu_t}
\newcommand{\Hmuinfty}{{\bf H}\mu_\infty}

\newcommand{\F}{{\mathscr F}}
\newcommand{\FF}{{\mathscr F}C_b^2(\D(A\s))}

\begin{document}

\baselineskip 13pt

\title[Transition semigroups of O-U processes]
{Transition Semigroups of Banach Space Valued Ornstein-Uhlenbeck Processes$^1$}

\author{B. Goldys}
\address{School of Mathematics\\
The University of New South Wales\\
Sydney 2052, Australia}
\email{B.Goldys@unsw.edu.au}

\author{J.M.A.M. van Neerven} 
\address{Department of Applied Mathematical Analysis\\
Technical University of Delft \\ P.O. Box 5031\\ 2600 GA Delft\\The Netherlands}
\email{J.M.A.M.vanNeerven@tudelft.nl}




\begin{abstract} We investigate the 
transition semigroup of the solution to a stochastic evolution
equation 
$$ dX(t) = AX(t)\,dt +dW_H(t),\qquad t\ge 0,$$ 
where $A$ is the generator of a $C_0$-semigroup $\SS$
on a separable real Banach space $E$ and
$\{W_H(t)\}_{t\ge 0}$ is cylindrical white noise with values in
a real Hilbert space $H$ which is continuously
embedded in $E$. Various properties of these semigroups,
such as the strong Feller property, the spectral gap property, and analyticity,
are characterized in terms of the behaviour of $\SS$ in $H$. In particular we 
investigate the interplay between analyticity of the transition semigroup,
$\SS$-invariance of $H$, and analyticity of the restricted 
semigroup $\SS_H$. 
\end{abstract}

\maketitle
\tableofcontents

\setcounter{footnote}{1}
\footnotetext{This is an update of the published version in Acta Appl. Math. {\bf 76}  (2003), 283--330. Changes are marked with footnoted}

\section{Introduction}
     In this paper we study transition semigroups associated
with stochastic linear Cauchy problems
\begin{equation}\label{bg0}
\begin{aligned}
dX(t) & = AX(t)+dW_H(t),\qquad t\ge 0,\\
 X(0) & = x.\end{aligned}
\end{equation}
We assume that $A$ is the generator of a $C_0$-semigroup
${\bf S}=\{S(t)\}_{t\ge 0}$ of bounded linear operators on a separable
real Banach space $E$ and $\WW_H = \left\{W_H(t)\right\}_{t\ge 0}$ is a cylindrical
Wiener process with the Cameron-Martin space $H$ which
is continuously imbedded into $E$. 

If $E$ is a Hilbert space, 
an explicit condition is known (see for example
\cite{DZ}) which ensures the existence of a unique
solution to (\ref{bg0}) of the form
\begin{equation}X(t,x)=S(t)x+\int_0^tS(t-s)\,dW_H(s).\label{ii}\end{equation}
The solution $\left\{X(t,x)\right\}_{t\ge 0}$ is 
called the {\em Ornstein-Uhlenbeck \/}
process associated with $\SS$ and $\WW_H$. It is a Markov process on $E$ 
whose transition semigroup is given by
\begin{equation}
P(t)\phi (x)=\mathbb E\phi (X(t,x))=\int_E\phi\left(S(t)x+y\right)\,d\mu_t(y),
\label{iii}
\end{equation}
where $\left\{\mu_t:t\ge 0 \right\}$ is a family of centred Gaussian
measures on $E$ associated with $\SS$ and $H$;
see Section \ref{sec:OU1} for details. 
This semigroup is also called the {\em Ornstein-Uhlenbeck semigroup}
associated with $\SS$ and $H$. 

If $E$ is a
Banach space, there seems to be no general satisfactory
theory of stochastic integration to give a rigorous meaning to the 
integral appearing in (\ref{ii}).
However, in many important cases
it can be shown that formula (\ref{ii}) is meaningful (at
least in a weak sense) and defines again a Markov
process on $E$ with transition semigroup ${\bf P}=\left\{P(t)\right
\}_{t\ge 0}$ given by (\ref{iii}); 
see for example \cite{BRS}, \cite{Brz}, \cite{BN}. The aim of this paper is to
study the transition semigroup ${\bf P}$ and its generator
under the sole assumption that the process 
$\left\{X(t,x)\right\}_{t\ge 0}$ is well
defined and admits an invariant measure $\mu_{\infty}$.

Apart from the case where $E$ is itself a Hilbert space and $H=E$,
many aspects of Ornstein-Uhlenbeck semigroups are not well understood.
For example, the existing criteria for the strong Feller property are difficult
to check in general. Similarly,
it is very difficult to check whether $\PP$ is analytic in
$L^2(E,\mu_\infty)$
or whether its generator has the spectral gap property.

The main idea of this paper, already exploited in \cite{CG2, CG3},
is to study the transition semigroup under the assumption that 
$\SS$ restricts to a $C_0$-semigroup $\SS_H$ on $H$. In this setting we obtain
explicit conditions for some properties of $\PP$ in terms of the behaviour of
the semigroup $\SS_H$. In particular we provide necessary and sufficient
conditions for the strong Feller property of $\PP$ and for the existence of 
a spectral gap. We also obtain conditions for analyticity of $\PP$
in terms of analyticity of the restricted semigroup $\SS_H$ which seem to be
close to optimal.

Our results extend and complement various
results from \cite{CG2, CG3, DP, DP2, DZ, fu, go, clos, Ne}.
 
\medskip
Let us now describe the contents of the paper in more detail.
Since many properties of ${\bf P}$ are determined by the
behaviour of the semigroup ${\bf S}$ on the spaces $H$ and
the reproducing kernel Hilbert spaces $H_t$ associated with the measures 
$\mu_t$, Sections \ref{sec:preliminaries} 
and \ref{sec:H-invariant}
are devoted to a study of interactions between the semigroup ${\bf S}$, the
space $H$ and the spaces $H_t$. We also investigate in
deatil the situation when $H$ is invariant under the
semigroup ${\bf S}$. In Section \ref{sec:lyapunov} the Liapunov
equation is considered and conditions are given for the
symmetry of ${\bf S}$ acting in $H$.

In Section \ref{sec:gap} we give several characterizations
of the spectral gap property of the generator $A$ of ${\bf S}$ when
considered in $H$ and reproducing kernel Hilbert space
$H_{\infty}$ associated with the invariant measure $\mu_\infty$.
In the case when $E$ is a Hilbert space,
it was shown in \cite{CG2} that this property is equivalent to the 
logarithmic Sobolev inequality
for the generator of the associated Ornstein-Uhlenbeck semigroup.
We also show that more accurate information can be obtained
if $H$ is $\SS$-invariant.

In Section \ref{sec:OU1} we introduce the Ornstein-Uhlenbeck
semigroup ${\bf P}$. It is studied in the space $C_b(E)$ endowed
with the mixed topology $\tau_{{\rm mixed}}$ under the minimal assumption
that (\ref{iii}) is meaningful. We extend the results from
\cite{GK} by showing that ${\bf P}$ is
$C_0$-semigroup in $\left(C_b(E),\tau_{{\rm mixed}}\right)$ and 
by giving an explicit formula for its generator $L$ on a suitable core. 
Let us note that ${\bf P}$ is not
strongly continuous, if fact not even strongly measurable, in $C_b(E)$
endowed with the supremum norm. We also provide a new
explicit condition for the strong Feller property of ${\bf P}$ in
the case when $H$ is ${\bf S}$-invariant.

Under the assumption of the existence of an invariant measure $\mu_\infty$, 
in Section \ref{sec:OU2} we study the semigroup ${\bf P}$ in $L^2\left(E,\mu_{
\infty}\right)$. In
particular we extend the existing criteria for the
symmetry of ${\bf P}$ and the existence of spectral gap for $L$.

In Sections
\ref{sec:analytic} and \ref{sec:analytic_H1} we are concerned with analyticity of
the Ornstein-Uhlenbeck semigroup in $L^2(E,\mu_\infty)$. 
We obtain necessary and sufficient conditions
for analyticity in terms of $H$. 
We establish connections between the analyticity of $\PP$,
the invariance of $H$ under $\SS$, and the analyticity of the 
restricted semigroup $\SS_H$.
We apply our criteria to prove analyticity of Ornstein-Uhlenbeck semigroups
associated with some stochastic partial differential equations of parabolic type.

\section{Preliminaries}
\label{sec:preliminaries}

\subsection{Reproducing kernel Hilbert spaces}

Stochastic evolution equations in Banach spaces are studied
conveniently by using the language of reproducing kernel Hilbert spaces.
We start by recalling some elementary properties of these spaces.

Throughout this paper $E$ denotes 
a real Banach space. The dual of $E$ is denoted by
$E\s$. A bounded linear operator 
$Q\in \L(E\s,E\s)$ is called {\em positive} if
$$
\lb Qx\s, x\s\rb \ge 0,\qquad x\s\in E\s,
$$
 and {\em symmetric} if
$$
 \lb Qx\s, y\s\rb = \lb Qy\s,x\s\rb, \qquad x\s, y\s\in E\s.
$$
More generally these definitions make sense for operators $Q\in\L(E\s,E\s\s)$.

If $Q\in\L(E\s,E)$ is positive and symmetric, then the bilinear map
on the range of $Q$ defined by
$$(Qx\s, Qy\s) \mapsto \lb Qx\s, y\s\rb, \qquad x\s,y\s\in E\s,$$
is easily checked to be a well defined inner product on the range of $Q$.
The Hilbert space completion of range$\,Q$ with respect to 
this inner product is called the {\em reproducing kernel Hilbert space} (RKHS) associated with $Q$ and is
denoted by $(H_Q, [\cdot,\cdot]_{H_Q})$. 
It is well known that 
the inclusion mapping from range$\,Q$ into $E$ extends to a
continuous injection from $H_Q$ into $E$. Denoting this extension by $i_Q$, we
have
$$Q = i_Q\circ i_Q\s.$$
This factorization immediately implies that $Q$ 
is weak$\s$-to-weakly continuous, and that  $H_Q$ is separable
whenever $E$ is separable.

Conversely, if $i:H\embed E$ is a continuous embedding of a real Hilbert space 
$H$ into $E$, then $Q:= i\circ i\s$ is positive and symmetric. As subsets of
$E$ we have $H = H_Q$ and
the map $i\s x\s\mapsto i_Q\s x\s$ defines
an isometrical isomorphism of $H$ onto $H_Q$. 

\begin{example} \label{ex:RKHS}

\

\begin{enumerate}
\item If $B$ is a bounded operator from a real Hilbert space
${\mathscr H}$ into $E$, then $Q:=B\circ B\s\in\L(E\s,E)$ 
is positive and symmetric. 
As subsets of $E$ we have $H_Q = \hbox{range}\,B$, and the inner product 
of $H_Q$ is given by
$$[Bg, Bh]_{H_{Q}} = [Pg, Ph]_{\mathscr H},
\qquad g,h\in {\mathscr H}.$$
Here $P$ denotes the orthogonal projection in ${\mathscr H}$ onto the 
orthogonal complement of ker$\,B$.
\item
As a special case of (1) let $E$ be a real Hilbert space
and let $Q\in\L(E)$ be a positive and selfadjoint operator.
Identifying the dual space $E\s$ with $E$ in the natural
way, we have
$H_Q = \hbox{range}\ Q^{\frac12}$, with inner product
$$[Q^{\frac12}x, Q^{\frac12}y]_{H_Q} = [Px,Py]_E, \qquad x,y\in E.$$
Here $P$ denotes the orthogonal projection in $E$ onto the 
orthogonal complement of ker$\,Q^\frac12$.
\end{enumerate}
\end{example}

It will be useful to compare the RKHS's associated with 
different positive symmetric operators in $\L(E\s,E)$. In this direction
we have the following easy fact; cf. \cite[Appendix B]{DZ}.
If $Q$ and $R$ are positive and symmetric operators in $\L(E\s,E)$, 
the following assertions are equivalent:
\begin{enumerate}
\item $i_Q(H_Q)\subseteq i_R(H_R)$;
\item There exists a constant $M\ge 0$ such that 
\begin{equation}\label{eq:incl}
\lb Qx\s,x\s\rb \le M \lb Rx\s,x\s\rb, \qquad x\s\in E\s.
\end{equation}
\end{enumerate}
Whenever it is convenient, we shall identify an embedded Hilbert space
with its image in $E$. Thus, 
instead of $i_Q(H_Q)\subseteq i_R(H_R)$ we shall simply write 
$H_Q \subseteq H_{R}$.

Another simple observation about RKHS's will be useful. Suppose $E$ and $F$
are   real Banach spaces, $j:E\embed F$ a continuous inclusion, and 
$Q_E\in\L(E\s,E)$ and $Q_F\in\L(F\s,F)$ are positive symmetric operators 
such that the following diagram commutes:
\begin{equation*}
 \begin{CD}
    E     @>> j >   F  \\
     @AA Q_E A                @AA Q_F A \\
     E\s     @<< j\s <  F\s 
     \end{CD}
\end{equation*}
Thus,  $Q_F = j\circ Q_E\circ j\s$. 
Let $i_E: H_E \embed E$ and  $i_F: H_F \embed F$ denote the RKHS's 
associated with $Q_E$ and $Q_F$, respectively. 
Then the mapping
$$I_{E,F}: i_E\s j\s y\s \mapsto i_F\s y\s, \qquad y\s\in F\s,$$
extends uniquely to an isometry from $H_E$ onto $H_F$. Moreover,
as subsets of $F$, the spaces $H_E$ and $H_F$ are identical.

Indeed, we compute: 
$$\n  i_E\s j\s y\s \n_{H_E}^2 = \lb Q_E j\s y\s, j\s y\s\rb
= \lb Q_F y\s, y\s\rb = \n i_F\s y\s\n_{H_R}^2, \qquad y\s\in F\s.
$$ 
Since $i_E$ and $j$ are injective, $ i_E\s \circ j\s$ has dense range in $E\s$,
and since $i_F$ is injective, $ i_F\s$ has dense range in $F\s$.
This shows that $I_{E,F}$ uniquely extends to an isometry of $H_E$ onto $H_F$.
From $I_{E,F}\circ i_E\s\circ j\s = i_F\s$
it follows, moreover, that $j\circ i_E\circ I_{E,F}\s = i_F$
and therefore, 
$$j(i_E(H_E)) =(j\circ i_E \circ I_{E,F}\s \circ I_{E,F}\s )(H_E)
=  (i_F \circ I_{E,F}) (H_E) = i_F (H_F).$$
This shows that $H_E$ and $H_F$ are identical as subsets of $F$
and we obtain the commuting diagram

\begin{equation*}
 \begin{CD}
    E     @>> j >   F  \\
     @AA i_E A                @AA i_F A \\
    H_E     @>> \simeq >  H_F 
     \end{CD}
\end{equation*}

The following observation will be useful:

\begin{proposition}\label{E-valued}
Let $R\in \L(E\s, E\s\s)$ be a positive symmetric operator. 
Suppose there 
exists a positive symmetric operator $Q\in\L(E\s,E)$ and a constant $C\ge 0$ 
such that
$$ \lb x\s, Rx\s\rb \le C\, \lb Qx\s,x\s\rb, \qquad x\s\in E\s.$$
Then $R\in\L(E\s,E)$.
\end{proposition}
\begin{proof}
Fix $x\s\in E\s$.
On the range of $i_Q\s$ we define a linear form $\phi_{x\s}$ by
$$\phi_{x\s} ( i_Q\s y\s) := \lb y\s, R x\s\rb, \qquad y\s\in E\s.$$
By the Cauchy-Schwarz inequality applied to the symmetric bilinear form
$(x\s,y\s)\mapsto \lb y\s, Rx\s\rb$, 
$$
\begin{aligned}
| \phi_{x\s} ( i_Q\s y\s) | & = | \lb y\s, R x\s\rb| 
\le \lb  x\s, Rx\s\rb^\frac12 
       \lb  y\s,R y\s\rb^\frac12
\\ & \le 
C\lb Qx\s,x\s\rb^\frac12 \lb Qy\s,y\s\rb^\frac12 
  = \n i_Q x\s\n_{H_Q} \n i_Q\s y\s\n_{H_Q}.
\end{aligned}
$$
It follows that $\phi_{x\s}$ is well defined and 
extends to a bounded linear form on $H_Q$
of norm $\le{C}\n i_Q\s x\s\n_{H_Q}$. By the Riesz representation
theorem, we may identify $\phi_{x\s}$ with an element of $H_Q$.
For all $y\s\in E\s$ we then have
$$ \lb i_Q \phi_{x\s}, y\s\rb = [\phi_{x\s}, i_Q\s y\s]_{H_Q}
= \phi_{x\s} (i_Q\s y\s) = \lb y\s, R x\s\rb.
$$
This shows that  $ Rx\s = i_Q \phi_{x\s} \in E$.
\end{proof}

\subsection{The general setting}

We consider a $C_0$-semigroup 
$\SS = \{S(t)\}_{t\ge 0}$ of boun\-ded linear operators on $E$
and a   real Hilbert space $H$ which is continuously embedded into $E$.
The embedding will be denoted by $i:H\embed E$. The inner product of $H$
will be denoted by $[\cdot,\cdot]_H$. 
The operator 
$Q := i\circ  i\s \in \L(E\s,E)$ is positive and symmetric, and $H$ is its
RKHS.

By \cite[Proposition 1.2]{Ne},
the $E$-valued function
$s\mapsto S(s)QS\s(s)x\s$ is strong\-ly measurable and we may define, for each 
$t>0$, the positive symmetric operator
$Q_t\in \L(E\s,E)$ by
\begin{equation}
\label{eq:defQ_t}
 Q_t x\s := \int^t_0 S(s)QS\s(s)x\s\,ds, \qquad x\s\in E\s.
\end{equation}
The RKHS associated with $Q_t$ will be denoted by $H_t$ and the embedding
$H_t\embed E$ by $i_t$. 
From \eqref{eq:incl}
it is immediate that
$H_s\subseteq H_t$ whenever $s\le t$ and 
the inclusion mapping is contractive \cite[Corollary 1.5]{Ne}.
Whenever it is convenient we further put $Q_0 := 0$ and $H_0 = \{0\}$.

\medskip

We will frequently consider the following hypothesis:
\begin{itemize}

\item $(\HQinfty)$: For all $x\s\in E\s$, {\em
weak}\,-\,$\lim_{t\to\infty} Q_t x\s $
exists in $E$.
\smallskip

\end{itemize}
Here, `{\em weak}\,-\,$\lim$' denotes the limit in the weak topology of $E$. 
This hypothesis is slightly more general than the one in 
\cite[Section 6]{Ne} where strong limits are taken, but the results 
proved there remain true under ($\HQinfty$)
without change in the proofs.

Assuming  ($\HQinfty$), we may define a bounded
operator $\Q: E\s\to E$  by
$$\Q x\s := weak\hbox{-}\lim_{t\to\infty} Q_t x\s, \qquad x\s\in E\s.$$ 
Clearly, $\Q$ is positive and symmetric. 
The RKHS associated with $\Q$ will be denoted by $\H$ and the embedding
$\H\embed E$ by $i_\infty$. 
From \eqref{eq:incl}
it is immediate that 
$H_t\subseteq \H$ for all $t>0$; by an obvious modification of
\cite[Corollary 1.5]{Ne} the inclusion mapping is contractive.

Necessary and sufficient conditions for $(\HQinfty)$ to be satisfied 
will be given in Section \ref{sec:lyapunov}. 
Hypothesis $(\HQinfty)$ is
trivially satisfied if $\SS$ is {\em uniformly exponentially stable},
i.e. if there exist constants $M\ge 0$ and $\o>0$ such that
$\n S(t)\n \le Me^{-\o t}$ for all $ t\ge 0.$
In this case we have 
\begin{equation}
\label{bochner}
\Q x\s = \int^\infty_0 S(s)QS\s(s)x\s\,ds, \qquad x\s\in
E\s,
\end{equation}
the integral being convergent as a Bochner integral in $E$.

Even in the case when $E$ is separable,  
we do not know whether the integral in \eqref{bochner}
always exists as a Bochner integral. 
We shall prove next that the integral always 
does exist as a Pettis integral.
For more information on Pettis integrals we refer the reader to \cite{DU}.

\begin{proposition}[$\HQinfty$]
For all $x\s\in E\s$ we have
%
$$
\Q x\s = \int^\infty_0 S(s)QS\s(s)x\s\,ds,
$$
the integral being convergent as a Pettis integral in $E$.
\end{proposition}
\proof 
Let $x\s\in E\s$ be fixed. First we prove the following claim: for all 
$y\s\in E\s$ the real-valued 
function $s\mapsto \lb
S(s)QS\s(s)x\s,y\s\rb$ is Lebesgue integrable on $[0,\infty)$ 
and 
\begin{equation}
\label{leb}
\lb Q_\infty x\s, y\s\rb = \int^\infty_0\lb S(s)QS\s(s)x\s,y\s\rb\,ds.
\end{equation}
First we take $y\s = x\s$. We have
$$ \lb S(s)QS\s(s)x\s,x\s\rb = \lb QS\s(s)x\s, S\s(s)x\s\rb \ge 0,
\qquad s\ge 0.$$
Hence by monotone convergence,
$$
\begin{aligned}
& \int^\infty_0|\lb S(s)QS\s(s)x\s,x\s\rb|\,ds
 = \int^\infty_0\lb S(s)QS\s(s)x\s,x\s\rb\,ds
\\ & \qquad \qquad = \lim_{t\to\infty} \int^t_0 \lb S(s)QS\s(s)x\s,x\s\rb\,ds
  = \lim_{t\to\infty} \lb Q_t x\s,x\s\rb
= \lb Q_\infty x\s,x\s\rb.
\end{aligned}
$$
This proves the claim for  $y\s = x\s$. 

Next let $y\s\in E\s$ be arbitrary. For all $t >0 $ we have
$$
\begin{aligned}
 & \int^t_0 |\lb  S(s) QS\s(s)x\s, y\s\rb|\,ds 
 = \int^t_0 |[ i\s S\s(s)x\s, i\s S\s(s)y\s]_H|\,ds
\\ & \qquad \qquad \le \left( \int^t_0 \n i\s S\s(s)x\s\n_H^2\,ds\right)^\frac12  
 \cdot \left(\int^t_0 \n i\s S\s(s)y\s\n_H^2\,ds \right)^\frac12 
\\ & \qquad  \qquad= \left( \int^t_0 \lb Q S\s(s)x\s, S\s(s)x\s\rb\,ds\right)^\frac12  
 \cdot \left(\int^t_0 \lb Q S\s(s)y\s, S\s(s)y\s\rb \,ds\right)^\frac12 
\\ & \qquad  \qquad\le \left( \int^\infty_0 \lb Q S\s(s)x\s,
S\s(s)x\s\rb\,ds \right)^\frac12  
 \cdot \left(\int^\infty_0 \lb Q S\s(s)y\s, S\s(s)y\s\rb\,ds \right)^\frac12 
\\ & \qquad \qquad =  \lb Q_\infty x\s,x\s\rb^\frac12\cdot 
\lb Q_\infty y\s,y\s\rb^\frac12.
\end{aligned}
$$
Passing to the limit $t\to\infty$, we obtain
$$  \int^\infty_0 |\lb S(s)QS\s(s)x\s, y\s\rb|\,ds 
\le \lb Q_\infty x\s,x\s\rb^\frac12\cdot 
\lb Q_\infty y\s,y\s\rb^\frac12.
$$
It follows that $s\mapsto \lb S(s)QS\s(s)x\s, y\s\rb$ is 
Lebesgue integrable in $[0,\infty)$.
The identity \eqref{leb} now follows from the dominated convergence theorem.
This concludes the proof of the claim.

In order to prove that $t\mapsto S(t)QS\s(t)x\s$ is Pettis integrable,
we have to show next that for all measurable subsets $B\subseteq [0,\infty)$
there exists an element $x_{B,x\s}\in E$ such that 
$$ \lb x_{B,x\s}, y\s\rb = \int_B \lb S(t)QS\s(t)x\s, y\s\rb\,dt, \qquad y\s\in E\s.$$
To this end, define the positive symmetric operator 
$Q_B\in\L(E\s,E\s\s)$ by
$$ \lb y\s, Q_B x\s\rb := \int_B \lb S(t)QS\s(t)x\s, y\s\rb\,dt.\qquad
x\s,y\s\in E\s.$$
Clearly, for all $x\s\in E\s$ we have $\lb Q_B x\s, x\s\rb \le \lb \Q
x\s,x\s\rb$, and therefore $Q_B\in\L(E\s,E)$ by Proposition
\ref{E-valued}. Then $x_{B,x\s}:= Q_B x\s$ does the job.
\qed

The space $\H$ displays some remarkable properties, some of which we shall
discuss next.

\begin{proposition}[$\HQinfty$]\label{prop:Hinfty-1}
The space $\H$ is invariant under the action of $\SS$, and the
restriction $\SS_\infty$ of $\SS$ to $\H$ defines a strongly continuous
contraction semigroup on $\H$. Its adjoint $\SS_\infty\s$ is {\em strongly
stable}, i.e.  for all $h_\infty\in\H$ we have
$$\lim_{t\to\infty} \n S_\infty\s(t)h_\infty\n_{\H} = 0.$$
\end{proposition}
\proof
The first assertion is proved in \cite{CG1} (for Hilbert spaces $E$) and 
\cite{Ne}. 

Noting that $S(t)\circ \i = \i \circ S_\infty(t)$, for all $x\s\in E\s$ we have
$$
\begin{aligned}
 \lim_{t\to\infty} \n S_\infty\s (t) i_\infty\s x\s\n_{\H}^2
 & =\lim_{t\to\infty} \n i_\infty\s S\s(t) x\s\n_{\H}^2
= \lim_{t\to\infty} \lb \Q S\s(t) x\s, S\s(t)x\s\rb
\\  & =\lim_{t\to\infty} \int^\infty_0 \lb S(s)Q
S\s(s)S\s(t)x\s,S\s(t) x\s\rb \,ds
\\ & = \lim_{t\to\infty}  
    \int^\infty_t \lb S(\sigma)Q S\s(\sigma)x\s,x\s\rb \,d\sigma
=0.
\end{aligned}
$$
Since the range of $i_\infty\s$ is dense in $\H$ and 
$\SS_\infty$ is a contraction semigroup on $\H$, the strong stability
of $\SS_\infty\s$ follows from this.
\qed

For later reference we recall from \cite{CG1} and \cite{Ne}:
 
\begin{proposition}[$\HQinfty$]\label{prop:Hinfty-2}
For $t>0$ fixed, the following assertions are equivalent:
\begin{enumerate}
\item $H_t = H_\infty$ with equivalent norms;
\item $\n S_\infty(t)\n_{\H} <1$.
\end{enumerate}
\end{proposition}

The following result gives a relation between $\H$ and $H$:

\begin{proposition}[$\HQinfty$]
We have $\overline{H}\,\subseteq\, \overline{\H}$, the
closures being taken in $E$. 
\end{proposition}
\proof
Suppose $y\s\in E\s$ is such that $\lb h_\infty ,y\s\rb = 0$ for all
$h_\infty \in\H$; we have to
prove that $\lb h,y\s\rb = 0$ for all $h\in H$.

First note that 
from $H_t\subseteq H_\infty$ it follows that $\lb Q_t x\s,y\s\rb = 0$
for all $t>0$ and $x\s\in E\s$.
Now fix $x\z\in E\z$, where $E\z$ denotes 
the closed linear subspace of $E\s$
of all elements whose orbit under the adjoint semigroup $\SS\s$ is
strongly continuous. Then for all $t>0$ we have
$$ \int^t_0 \lb S(s)QS\s(s)x\z,y\s\rb\,ds =\lb Q_t x\z, y\s\rb =0,$$ 
and since the integrand is a continuous function, this implies that
$$\lb S(s)QS\s(s)x\z,y\s\rb = 0, \qquad s\ge 0.$$ 
In particular, $\lb Qx\z, y\s\rb = 0$.

Since $Q$ is symmetric, it follows that $\lb Qy\s, x\z\rb = 0$
for all $x\z\in E\z$, and $E\z$ being weak$\s$-dense in $E\s$ this
implies that $Qy\s=0$. Then $\lb Qx\s, y\s\rb = \lb Qy\s,x\s\rb=0$ for all 
$x\s\in E\s$, and since the range of $Q$ is dense in $H$
it follows that  $\lb h,y\s\rb = 0$ for all $h\in H$.
\qed

It need not be the case that $H\subseteq \H$. In fact,
as we will show in Section \ref{sec:gap} it often happens that
$\H\subseteq H$ (in which case of course $\overline{H} = \overline{\H}$).

\section{Invariance of the reproducing kernel Hilbert space $H$}
\label{sec:H-invariant}

In many important examples, $H$ is invariant under the action of $\SS$ and
$\SS$ restricts to a $C_0$-semigroup on $H$. For example, 
we will show that this happens if the 
Ornstein-Uhlenbeck semigroup $\PP$ in $L^2(E,\mu_\infty)$ 
is selfadjoint (Section \ref{sec:lyapunov}) or analytic
with a spectral gap (Section \ref{sec:analytic_H1}).
A further example is when $E$ is a Hilbert space and $S(t)Q=QS(t)$ holds 
for all $t\ge 0$; see \cite{CG2}. 

In this section we will investigate the situation where $\SS$ restricts to a
$C_0-$semi\-group on $H$ in some detail. It will turn out that the 
restricted semigroup enjoys some
interesting regularizing properties. These will be used to study the 
strong Feller property of Ornstein-Uhlenbeck semigroups.

\medskip
We begin with a simple criterion for invariance.
If $T\in\L(E)$ is a bounded operator satisfying $T(H)\subseteq H$,
then we denote the
restriction of $T$ to $H$ by $T_H$; by the closed graph theorem, 
$T_H$  is a bounded operator on $H$.
Note that $T\circ i = i\circ T_H$.

\begin{proposition}\label{Riesz}
For a bounded operator  $T\in\L(E)$ the following assertions are equivalent:

\begin{enumerate}
\item $T(H)\subseteq H$;
\item There exists a constant $M\ge 0$ such that for all $x\s\in E\s$
we have
$$ \n i\s T\s x\s\n_H \le M \n i\s x\s\n_H.$$ 
\item There exists a constant $M\ge 0$ such that for all $x\s, y\s\in E\s$ 
we have $$|\lb TQx\s, y\s\rb| \le M \n i\s x\s\n_H \n i\s y\s\n_H.$$
\end{enumerate}
In this situation the restriction $T_H$
is bounded on $H$ and satisfies $\n T_H\n_H\le M$, where $M$ is either one
of the constants in $(2)$ or $(3)$.
\end{proposition}


\begin{proof}
(1)$\,\Rightarrow\,$(2):
From $T\circ i = i\circ T_H$ we have, for all $x\s\in E\s$,
$$\n i\s T\s x\s\n_H = \n T_H\s (i\s x\s)\n_H \le \n T_H\n_H \n i\s x\s\n_H.$$
This gives (2), with $M = \n T_H\n_H$.

(2)$\,\Rightarrow\,$(3):
From $Q = i\circ i\s$ we then have,  for all $x\s,y\s\in E\s$,
$$ |\lb TQx\s, y\s\rb| = |\lb i (i\s x\s), T\s y\s\rb| = |[i\s x\s, i\s T\s y\s]_H|
\le \n T_H\n_H  \n i\s x\s\n_H \n i\s y\s\n_H.$$
This gives (3), with the same constant $M$.

(3)$\,\Rightarrow\,$(1):
By assumption, the mapping 
$\phi: i\s y\s \mapsto \lb TQx\s, y\s\rb$ is well defined and 
uniquely extends to a bounded 
linear functional
$\phi$ on $H$ of norm $\le M\n i\s x\s\n_H$. 
By the Riesz representation theorem we identify $\phi$
with an element $h\in H$ of norm  $\le M\n i\s x\s\n_H$. 
Then for all $y\s\in E\s$ we have
$$ \lb ih, y\s \rb = [i\s y\s, h]_H  =  \phi(i\s y\s) = \lb TQx\s, y\s\rb,$$
and therefore $TQ x\s = ih.$
Defining $T_H(i\s x\s) := h$, we have $\n T_H(i\s x\s)\n = \n h\n_H
\le  M\n i\s x\s\n_H$. Hence we obtain a well defined
bounded operator $T_H$ on $H$ of norm
$\le M$.
Finally, for all $x\s,y\s\in E\s$ we have
$$ \lb (i\circ T_H)(i\s x\s), y\s\rb  = \lb ih, y\s\rb = \lb TQx\s, y\s\rb = 
\lb(T\circ i)(i\s x\s),y\s\rb,
$$
which shows that $i\circ T_H = T\circ i$.
\end{proof}

The implication (2)$\,\Rightarrow\,$(1) admits the following,
even shorter, direct proof.
By assumption of (2), the mapping $S_H: i\s x\s \mapsto  i\s T\s x\s$
is well defined and extends uniquely to a bounded operator
$S_H$ on $H$ of norm $\le M$.
From $S_H\circ i\s = i\s \circ T\s$ we obtain, by dualizing,
$i\circ S_H\s = T\circ i$.  
Hence $T$ maps $H$ into itself and 
the restriction of $T$ to $H$ equals $S_H\s$.
\medskip
Next we address the question of strong continuity of $\SS_H$.

\begin{proposition}\label{prop:str-cont}
Assume that $S(t)H\subseteq H$ for all $t\ge 0$.
Then the semigroup $\SS_H$ is strongly continuous on $[0,\infty)$ 
if and only if
$$\limsup_{t\downarrow 0} \n S_H(t)\n_H < \infty.$$
\end{proposition}
\proof
%
%
For all $x\s, y\s \in E\s$ we have
$$
 \begin{aligned}
\lim_{t\downarrow 0}\, [S_H(t) i\s x\s - i\s x\s, i\s y\s]_H 
& = \lim_{t\downarrow 0}\, [i\s x\s, i\s S\s(t) y\s - i\s y\s ]_H =0
 \\ &= \lim_{t\downarrow 0}\, \lb Qx\s, S\s(t) y\s - y\s \rb = 0.
\end{aligned}
$$
By assumption, $\SS_H$ is locally bounded near $0$, and it follows 
that $\SS_H$ is weakly continuous. By a standard result from semigroup 
theory, this implies that that $\SS_H$ is strongly continuous.
\qed

It may happen that $\SS_H$ fails to be
strongly continuous at $0$, even if $E$ is a Hilbert space:

\begin{example}\label{ex:notC0}
For $n = 1,2,\dots$ we define the Hilbert space $H_n$
to be $L^2[0,1]$ with the norm
$$ \n f\n_{H_n}^2 := \int^{1-\frac1n}_0 |f(t)|^2\,dt + n^2 
\int^{1-\frac1{2n}}_{1-\frac1n} |f(t)|^2\,dt
+  \int^1_{1-\frac1{2n}} |f(t)|^2\,dt.
$$
The nilpotent left shift semigroup $\SS_{H_n}$ on $H_n$,
 $$S_{H_n}(t)f(s) = 
\left\{
\begin{array}{cl}  
 f(s+t), & 0\le s+t\le 1 \\
   0   , & \hbox{else}
\end{array}
\right.
$$
is strongly continuous and we have
$$ 
\n S_{H_n}(t)f(s)\n_{H_n} = 
\left\{
\begin{array}{cl}  
   n   , & t\in [0,\tfrac1n) \\
   1   , & t\in [\tfrac1n, 1) \\
   0   , & \hbox{else.}
\end{array}
\right.
$$
Let $E_n = L^2[0,1]$ with the usual norm and
let $\SS_{E_n}$ denote the nilpotent left shift semigroup on $E_n$.
Now consider the Hilbert space direct sums
$$ H := \bigoplus_{n=1}^\infty H_n, \qquad E := \bigoplus_{n=1}^\infty E_n.$$
Note that $H\subseteq E$ with a continuous inclusion map, which we denote by $i$.

The semigroups $\SS_H := \bigoplus_{n=1}^\infty \SS_{H_n}$ and
$\SS_E := \bigoplus_{n=1}^\infty \SS_{E_n}$
act in $H$ and $E$, respectively, and $\SS_H$ is the restriction to $H$
of $\SS_E$. We have
$$
\n S_H(t)f(s)\n_{H} = 
\left\{
\begin{array}{cl}
   1   , & t=0  \\
   n   , & t\in [\frac{1}{n+1},\tfrac1n), \qquad n=1,2,\dots \\
   0   , & \hbox{else.}
\end{array}
\right.
$$
Thus, $\limsup_{t\downarrow 0} \n S_H(t)\n_H = \infty$ and $\SS_H$
fails to be strongly continuous in $H$ at $0$.
On the other hand $\SS_E$ is strongly continuous at $0$. 
\end{example}

The infinitesimal generator of $\SS$ will be denoted by $A$.
The next result gives necessary and sufficient conditions for
$\SS_H$ to be contractive:

\begin{theorem}\label{thm:contraction}
The following assertions are equivalent:

\begin{enumerate}
\item For all $t\ge 0$ we have $S(t)H\subseteq H$, and $\SS_H$ is a $C_0$-semigroup of contractions
on $H$;
\item For all $x\s\in E\s$, the function $t\mapsto \n i\s S\s(t)x\s\n_H$ is 
nonincreasing on $[0,\infty)$;
\item For all $x\s\in \D(A\s)$, the domain of $A\s$, 
we have $-\lb Qx\s, A\s x\s\rb \ge 0.$ 
\end{enumerate}
\end{theorem}

\proof
(1)$\,\Rightarrow\,$(3): 
Fix $x\s\in \D(A\s)$. Then for all $h\in H$,
$$\lim_{t\downarrow 0} \, \frac1t[S_H(t)i\s x\s - i\s x\s, h]_H
= \lim_{t\downarrow 0} \, \frac1t\lb ih, S\s(t) x\s - x\s\rb 
= \lb ih, A\s x\s\rb = [i\s A\s x\s, h]_H.
$$ Let $A_H$ denote the infinitesimal generator of $\SS_H$.
By a standard result from semigroup theory \cite[Theorem 2.1.3]{Pa}, 
the above identities imply that
$i\s x\s\in \D(A_H\s)$ and $A_H\s(i\s x\s) = i\s A\s
x\s$. The fact that $A_H$ generates a contraction semigroup on 
$H$ then gives, using \cite[Theorem 1.4.3]{Pa},
$$
 -\lb Qx\s, A\s x\s\rb = -[i\s x\s, i\s A\s x\s]_H 
 = -[i\s x\s, A_H\s (i\s x\s)]_H \ge 0.
$$

(3)$\,\Rightarrow\,$(2):
First consider a fixed $x\s\in \D(A\s)$. 
It is easy to see that the function $u\mapsto 
\n i\s S\s(u)x\s \n_H^2 = \lb QS\s(u)x\s, S\s(s)x\s\rb$
is differentiable and 
$$ 
\begin{aligned}
\frac{d}{du} \n i\s S\s(u)x\s \n_H^2 
& = 
\lb QS\s(u)A\s x\s, S\s(u)x\s\rb + \lb QS\s(u)x\s, S\s(u)A\s x\s\rb
\\ & = 2 \lb QS\s(u)x\s, A\s S\s(u)x\s\rb,
\end{aligned}
$$
where we used the symmetry of $Q$. Hence for all $t\ge s\ge 0$ we have
$$
\begin{aligned}
\n i\s S\s(s)x\s \n_H^2 - \n i\s S\s(t)x\s \n_H^2
& = -\int^t_s \frac{d}{du} \n i\s S\s(u)x\s\n_H^2\,du
\\ & = -2 \int^t_s \lb QS\s(u)x\s,A\s S\s(u)x\s\rb \,du\ge 0.
\end{aligned}
$$
For general $x\s\in E\s$, the result follows by approximation: noting that
$$ 
\n i\s S\s(\tau)y\s\n_H 
= \lim_{\l\to\infty} \l \n (\l-A_H\s)^{-1} i\s S\s(\tau)y\s\n_H
= \lim_{\l\to\infty} \l \n i\s S\s(\tau)(\l-A\s)^{-1}y\s\n_H,
$$
we can apply the above to $x\s = (\l-A\s)^{-1}y\s \in\D(A\s).$

(2)$\,\Rightarrow\,$(1):
Since $t\mapsto \n i\s S\s(t)x\s\n_H$ is nonincreasing on $[0,\infty)$,
for all $t\in [0,\infty)$ and all $x\s \in E\s$ we have 
$\n  i\s S\s(t)x\s\n_H\le \n i\s x\s\n_H$. 
Then from Proposition \ref{Riesz} it follows that the operators $S(t)$
restrict to contractions on $H$.
%
%
The strong continuity of $\SS_H$ follows from Proposition
\ref{prop:str-cont}.
\qed  

Recalling that $i:H\embed E$ denotes the inclusion map, we define the positive
symmetric operator $Q\in\L(E\s,E)$ by
$$ Q:= i\circ i\s.$$
Using this operator, we define the positive symmetric operators
$Q_t\in\L(E,E\s)$ by \eqref{eq:defQ_t}.
As before we let $H_t$ be the RKHS associated with $Q_t$, and
$i_t:H_t\embed E$ is the natural inclusion mapping. 

\begin{theorem}\label{thm:regularH}
 Assume that $\SS$ restricts to a $C_0$-semigroup 
$\SS_H$ on $H$.
\label{thm:H-invar}

\

\begin{enumerate}
\item For all $t,s>0$ we have $ H_t=H_s$ with equivalent norms;
\item For all $t>0$ we have $H_t\subseteq H$ with dense inclusion; 
\item For all $t>0$ we have $S(t)H\subseteq H_t$ and  
$$ \limsup_{t\downarrow 0}\, \sqrt{t}\, \n S(t)\n_{{\mathscr L}(H, H_t)} 
   \le \limsup_{t\downarrow 0}\n S_H(\cdot)\n_{H}.
$$
\end{enumerate}
\end{theorem}

\proof

We start with some general observations. 
For $t>0$
define the positive selfadjoint operator $R_t\in \L(H)$ by
$$
 R_t \, h := \int^t_0 S_H(s)S_H\s(s)h\,ds, \qquad h\in H.
$$
Let $G_t$ denote the RKHS associated with 
$R_t$ and let
$j_t: G_t\embed H$ denote the 
inclusion mapping. 
By \cite[Theorem 1.11]{Ne},
as subsets of $E$ we have
$H_t=G_t$. It follows that $H_t \subseteq H$. 

Denoting by $i_t:H_t\embed E$  the inclusion mapping, 
the map $$i_t\s x\s \mapsto j_t\s(i\s x\s), \qquad x\s\in E\s,$$
establishes an isometrical isomorphism of $H_t$ and $G_t$.
Thanks to this observation, 
in the rest of the proof we may identify $H_t$
with $G_t$.

For each $h\in H$, the function $f_{h}(s) := S_H(s)h$ 
belongs
to $L^2([0,t];H)$, and hence
\begin{equation}
\label{eq:control}
t\cdot S(t)h = \int^t_0 S(t)h\,ds 
= \int^t_0 S(t-s) f_h(s)\,ds.
\end{equation}
Hence by \cite[Appendix B]{DZ}, $S(t)h \in G_t$. 
This shows that $S(t)H\subseteq G_t = H_t$.

Let $$Y_t=\{S_H(s)h: \ 0<s\le t,\  h\in H\}.$$ 
By the strong continuity of $\SS_H$, 
the set $Y$ is dense in 
$H$. On the other hand, for all $0<s\le t$ we have
$H_s \subseteq H_t$ and therefore
$$S_H(s)h \in H_s \subseteq H_t.$$
It follows that $Y\subseteq
H_t\subseteq H$, and  $H_t$ is dense in $H$. This proves (2).

Fix $t_0>0$. From (2) we have
$$S(t_0)H_{t_0}\subseteq S(t_0)H \subseteq H_{t_0}.$$  Therefore 
by \cite[Theorem 1.9]{Ne}, for all $t\ge t_0$ we have 
$H_t = H_{t_0}$ with equivalent norms.
Since $t_0>0$ is arbitrary, this gives (1).

Fix $h\in H$ and $t_0>0$.
Using the language of control theory of
\cite[Appendix B]{DZ}, \eqref{eq:control} shows that 
the function $u(s) :=t_0^{-1} S_H(s)h$ 
is a control for reaching $S(t_0) h$ at time $t_0$.
The the minimum energy for a control to reach $S(t_0)h$ being equal 
$\n S(t_0)h\n_{H_{t_0}}$, it follows that 
$$
\n S(t_0)h\n_{H_{t_0}}^2
\le   \n u\n_{L^2([0,t_0];H)}^2 
= \frac{1}{t_0^2} \int^{t_0}_0 \n S_H(s)h\n_H^2\,ds.
$$
Therefore,
$$
\limsup_{t\downarrow 0}\, \sqrt t \, \n S(t)\n_{{\mathscr L}(H, H_t)} 
\le \limsup_{t\downarrow 0}\n S_H(\cdot)\n_{H}.
$$
This gives (3).
\qed

In case $E$ is a Hilbert space and $Q=I$, these estimates are
well known; cf. \cite{DZ}.

\medskip
Assertion (3) admits a control theoretic interpretation.
In order to explain this, we need to introduce some terminology.

Let ${\mathscr H}$ be a real Hilbert space,
let $B: {\mathscr H}\to {E} $ a bounded linear operator, and 
let $t_0>0$ be given.
We say that the pair $(\SS, B)$ is {\em null controllable
in time $t_0$} if for every $x\in {E}$
there exists a function $f\in L^2((0,t_0);{\mathscr H})$
such that the unique mild solution $u$ of the equation
$$
\begin{aligned}
u'(t) & = {A}u(t) + Bf(t), \qquad t\ge 0,\\
u(0)  & = x
\end{aligned}
$$
satisfies $$u(t_0)=0.$$ 
It is well known that the pair $({\SS }, B)$ is null controllable
in time $t_0$ if and only if ${S}(t_0)$ maps ${E}$ into $H_{R_{t_0}}$,
the RKHS associated with the positive symmetric operator
$R_{t_0} \in\L({E}\s, {E})$ defined by
\begin{equation}\label{eq:RT}
 R_{t_0} x\s := \int^{t_0}_0 {S}(s)BB\s{S}\s(s)x\s\,ds, 
                 \qquad x\s\in E\s,
\end{equation}
cf. \cite{DZ, Ne3}.

\begin{theorem}\label{thm:null} 
In the above situation, let there exist $\delta>0$ and a constant 
$M\ge 0$ such that
for all $t\in [0,\delta]$ and all $x\s\in E\s$ we have
\begin{equation}
\label{eq:null} 
\n B\s {S}\s(t)x\s\n_{\mathscr H} \le M \n B\s x\s\n_{\mathscr H}.
\end{equation}
Then the following assertions are equivalent:

\begin{enumerate}
\item The pair $({\SS},B)$ is null controllable for all $t>0$;
\item $S(t)E\subseteq\ ${\rm range}$\,B$ for all $t>0$.
\end{enumerate}
\end{theorem}

\begin{proof}
Let $R:= B\circ B\s$ and let $i_R:H_R\embed {E}$ denote the  RKHS associated
with $R$.  As outlined in Example \ref{ex:RKHS} (1), 
we may identify $H_R$ with the range of $B$. By Proposition \ref{Riesz}, 
the estimate \eqref{eq:null} 
implies that $H_R$ is $S(t)$-invariant for all 
$t\in [0,\delta]$, and that the restricted operators are
uniformly bounded on $H_R$. Then Proposition \ref{prop:str-cont} 
implies that ${\SS}$
restricts to a $C_0$-semigroup on $H_R$, and we are in a position to apply
Theorem \ref{thm:regularH}. 

(1)$\,\Rightarrow\,$(2): If   $({\SS}, B)$ 
is null controllable at time $t$,
then ${S}(t)(E) \subseteq H_{R_t} \subseteq H_R$, 
where the second inclusion
follows from Theorem \ref{thm:regularH} (2).

(2)$\,\Rightarrow\,$(1): For all $t>0$ we have 
${S}(t)x = {S}(\frac{t}{2})({S}(\frac{t}{2})x) 
\in {S}(\frac{t}{2})H_R 
\subseteq H_{R_\frac{t}{2}} = H_{R_t}$, where we used
Theorem \ref{thm:regularH} (1), (3).
\end{proof}

In the previous section we showed that Hypothesis ($\HQinfty$) implies
the inclusion
$\overline{H}\, \subseteq \, \overline{\H}$.
If $H$ is $\SS$-invariant we can prove more:

\begin{proposition}[$\HQinfty$]\label{prop:H-cap-Hinfty}
If $\SS$ restricts to a $C_0$-semigroup on $H$, then 
$H\cap \H$ is dense in both $H$ and $\H$.
In particular we have
$\overline{H}\, = \, \overline{\H}$.
\end{proposition}
\begin{proof}
By Theorem \ref{thm:regularH}, 
for all $h\in H$ we have $S_H(t)h\in H_t \subseteq H\cap \H$. 
Hence from $\lim_{t\downarrow 0} S_H(t) h = h$ strongly in $H$ it  
follows that $H\cap \H$ is dense in $H$.

For all $t>0$ and $x\s\in E\s$ we have $i_t\s x\s \in H_t \subseteq H\cap \H$.
We claim that $\lim_{t\to\infty}i_t\s x\s = \i x\s$ weakly in $\H$. 
Since the range of $\i\s$ is dense in $\H$, this will show that
$H\cap \H$ is weakly dense in $\H$, and therefore dense in $\H$.

To prove the claim we first recall from Section \ref{sec:preliminaries}
that the inclusion mapping
$H_t\embed \H$ is contractive. Therefore,
\begin{equation}
\label{bound_t}
 \n i_t\s x\s\n_{\H}^2 \le \n i_t\s x\s\n_{H_t}^2
= \lb Q_t x\s, x\s\rb \le \lb \Q x\s, x\s\rb = \n \i\s x\s\n_{\H}^2.
\end{equation}
Moreover, for all $y\s\in E\s$ we have
$$
\begin{aligned} 
\lim_{t\to\infty} [i_t\s x\s, \i\s y\s]_{\H} 
& = \lim_{t\to\infty} \lb \i i_t\s x\s, y\s\rb
 = \lim_{t\to\infty} \lb i_t i_t\s x\s, y\s\rb
\\ & = \lim_{t\to\infty} \lb Q_t x\s, y\s\rb
 = \lb \Q x\s, y\s\rb =   [\i\s x\s, \i\s y\s]_{\H}.
\end{aligned}
$$
Using once more the density of the range of $\i\s$, 
together with the uniform bound \eqref{bound_t}
this proves the claim.
\end{proof}

%

If $\SS_H$ is a $C_0$-semigroup of {\em normal} operators, 
then for individual orbits
we have the following version of the estimate in Theorem
\ref{thm:H-invar} (3):

\begin{theorem}\label{thm:normal}
If $\SS_H$ is a $C_0$-semigroup of normal operators on $H$, then
for all $t>0$ and $h\in H$ we have
$$
 \int^t_0 \n S_H(s) h\n_{H_t}^2\,ds = \n h\n_H^2.
$$
\end{theorem}
Note that the right hand side is independent of the semigroup $\SS$.

\proof
For $t>0$
define the positive selfadjoint operators $R_t, R_{\ast t}\in \L(H)$ by
$$
 R_t \, h := \int^t_0 S_H(s)S_H\s(s)h\,ds, \quad
 R_{\ast t} \, h := \int^t_0 S_H\s(s)S_H(s)h\,ds, \qquad h\in H.
$$
Let $G_t$ and $G_{\ast t}$ denote the RKHS's associated with 
$R_t$ and $ R_{\ast t}$, respectively, and let
$j_t: G_t\embed H$ and $j_{\ast t}: G_{\ast t}\embed H$ denote the 
inclusion mappings.
From
$$ \n j_t\s h\n_{G_t}^2 = \int^t_0  [S_H(s)S_H\s(s)h,h]_H\,ds
= \int^t_0 [S_H\s(s)S_H(s)h,h]_H\,ds = \n (j_{\ast t})\s h\n_{G_{\ast t}}^2
$$
it follows that $G_t$ and $G_{\ast t}$ are canonically isometrically isomorphic
as Hilbert spaces, and identical as subsets of $H$. 
Moreover, as we saw in the proof of Theorem \ref{thm:H-invar}, 
$H_t$ and $G_t$ are canonically isometrically isomorphic
as Hilbert spaces and we have
$H_t=G_t$ as subsets of $E$.
It follows that it suffices to prove that 
$$
 \int^t_0 \n S_H(s)h\n_{G_{\ast t}}^2\,ds = \n h\n_H^2.
$$

 From the normality of each operator $S_H(t)$ it is not difficult to
see that
\begin{equation}\label{eq:commute}
S_H(t) S_H\s(s) = S_H\s(s) S_H(t), \qquad t,s\ge 0.
\end{equation}
Indeed, by the semigroup property this is true whenever both $t$ and $s$ are 
integer multiples of a common fixed real number. The set of all pairs $(t,s)$
with this property being dense in $[0,\infty) \times [0,\infty)$, 
the general case follows by strong
continuity of $\SS_H$ and its adjoint $\SS_H\s$.

Using \eqref{eq:commute} we see that 
$$
R_{\ast t} S_H(\tau) = S_H(\tau) R_{\ast t}, \quad 
R_{\ast t} S_H\s(\tau) = S_H\s(\tau) R_{\ast t}, \qquad \tau\ge 0,
$$ 
and hence,  
\begin{equation}\label{eq:sqrt_commute}
 R_{\ast t}^{\frac12} S_H(\tau) = S_H(\tau) R_{\ast t}^{\frac12},\quad
 R_{\ast t}^{\frac12} S_H\s(\tau) = S_H\s(\tau) R_{\ast t}^{\frac12}, \qquad \tau\ge 0.
\end{equation}
Define the convolution operator $A_t$ from $L^2([0,t];H)$ into $H$ by
$$A_t \psi := \int^t_0 S_H\s(t-s)\psi(s)\,ds.$$
By \eqref{eq:sqrt_commute},
$$ R_{\ast t}^{\frac12}A_t = A_t R_{\ast t}^{\frac12}.$$
It is trivially checked that
$$ A_t\s h = S_H(t-\,\cdot)h, \qquad h\in H.$$
Hence, the kernel of $A_t$ is equal to the orthonormal complement in
$L^2([0,t];H)$ of the closed linear subspace $V_t$ 
spanned by all functions of the
form $s\mapsto S_H(t-s)h$, $h\in H$. Let $\pi_t$ denote the orthonormal
projection in $L^2([0,t];H)$ onto $V_t$.

Now fix $h\in H$. Then by Example \ref{ex:RKHS}(2), 
$R_{\ast t}^{\frac12}h\in G_{\ast t}, $ 
the RKHS associated with $R_{\ast t}$. 
Therefore by \cite[Appendix B]{DZ} there exists a function  
$\phi\in L^2([0,t];H)$ such that
$$
 R_{\ast t}^{\frac12}h = A_t \phi.
$$ 
Noting that $A_t \phi = A_t (\pi_t \phi)$, it follows that
\begin{equation}\label{eq:R_t-1}
 R_{\ast t} h = R_{\ast t}^{\frac12} A_t \phi =  R_{\ast t}^{\frac12}A_t (\pi_t\phi)
 = A_t (R_{\ast t}^{\frac12}(\pi_t\phi)).
\end{equation}
On the other hand,
\begin{equation}\label{eq:R_t-2}
 R_{\ast t} h = \int^t_0 S_H\s(t-s)S_H(t-s)h\,ds = A_t(S_H(t-\,\cdot)f).
 \end{equation}
Observing that $S_H(t-\,\cdot)h\in V_t$, and noting
that from $\pi_t \phi\in V_t$ it follows that also 
$R_{\ast t}^{\frac12}(\pi_t \phi) \in V_t$, from
\eqref{eq:R_t-1} and \eqref{eq:R_t-2} 
we now deduce that
$$
 S_H(t-\,\cdot)h = R_{\ast t}^{\frac12}(\pi_t\phi).
$$
But clearly, $R_{\ast t}^{\frac12}(\pi_t\phi) \in L^2([0,t];G_{\ast t})$.
It follows that $ S_H(\cdot)h\in L^2([0,t];G_{\ast t})$.

Finally, because $R_{\ast t}^{\frac12}$ is an isometry from $H$ onto $G_{\ast t}$
and $A_t$ is an isometry from $V_t$ onto $G_{\ast t}$,
$$
\begin{aligned}
 \int^t_0 \n S_H(s)h\n_{G_{\ast t}}^2\,ds
& = \int^t_0 \n S(t-s)h\n_{G_{\ast t}}^2\,ds 
\\ & = \int^t_0 \n R_{\ast t}^{\frac12}(\pi_t\phi(s))\n_{G_{\ast t}}^2 \,ds
 = \int^t_0 \n \pi_t\phi(s)\n_{H}^2 \,ds
\\ & = \n \pi_t\phi\n_{L^2([0,t];H)}^2
= \n A_t(\pi_t\phi)\n_{G_{\ast t}}^2
= \n R_{\ast t}^{\frac12}h\n_{G_{\ast t}}^2 
 = \n h\n_H^2.
\end{aligned}
$$
\qed

\begin{remark}
By \cite[Theorem 22.4.1]{HPh}, for normal semigroups we always have an estimate
\begin{equation}\label{bound_a}
\n S_H(t)\n_H \le e^{a t}, \qquad t\ge 0,
\end{equation}
for some $a\in \R$. 

Let us now assume that $\SS_H$ is an analytic semigroup which satisfies
\eqref{bound_a} for some $a<0$.
These assumptions imply that 
$ D_{A_H}(\hbox{$\frac12$},2) = \D((-A_H)^{\frac12}) $;
we use standard notations as can be found, e.g., in \cite{DZ}.
From this, in turn, it follows that there exists a constant $C\ge 0$ such that
\begin{equation}
\label{eq:D_A-inequality}
\int^t_0 \n S_H(s)h\n_{D_{A_H}(\frac12, 2)}^2\,ds\le C \n h\n_H, 
\qquad h\in H, \ \ t>0;
\end{equation}
see \cite[Appendix A]{DZ}.
On the other hand, for all $t>0$ we have
$H_t = D_{A_H}(\hbox{$\frac12$},2)$ with equivalent norms
see \cite[Appendix B]{DZ}. Therefore \eqref{eq:D_A-inequality} implies 
\begin{equation}\label{Ct}
\int^t_0 \n S_H(s)h\n_{H_t}^2\,ds\le C_t \n h\n_H, 
\qquad h\in H, \ \ t>0,
\end{equation}
with a constant $C_t$ depending on $t$. Theorem \ref{thm:normal}
shows that in the normal case
one has equality in \eqref{Ct} with $C_t = 1$. 
\end{remark}

\section{The Liapunov equation $AX + XA\s = -Q$ and $Q$-symmetry}
\label{sec:lyapunov}

In this section we study the {\em Liapunov equation}
\begin{equation}
\label{eq:lyapunov}
AX +XA\s = -Q
\end{equation}
and apply the results to the case where we have $S(t)\circ Q = Q\circ S\s(t)$
for all $t\ge 0$.

The following result shows that the operator $\Q$, if exists, 
`solves' this equation:

\begin{proposition}[$\HQinfty$]\label{prop:lyapun}
For all $x\s\in \D(A\s)$ we have $\Q x\s\in \D(A)$ and 
$$A\Q x\s + \Q A\s x\s = -Qx\s.$$
\end{proposition}
\proof
Take $x\s, y\s\in \D(A\s)$. 
Differentiating the identity
$$
 \lb \Q S\s(t)x\s, S\s(t) y\s\rb = \lb \Q x\s, y\s\rb - \lb Q_t x\s,y\s\rb
$$
on both sides with respect to $t$. Evaluating at $t=0$ gives
\begin{equation}\label{eq:lyapunov-Q-infty}
 \lb \Q x\s, A\s y\s\rb + \lb \Q A\s x\s, y\s\rb = - \lb Qx\s, y\s\rb.
\end{equation}
It follows that $\Q x\s\in \D(A)$ and that $A\Q x\s + QA\s x\s = -Q$.
\qed

This result motivates the following definition.

\begin{definition}
A {\em solution} of equation \eqref{eq:lyapunov} 
is a bounded operator $X\in \L(E\s,E)$ such that for 
all $x\s\in \D(A\s)$ we have
$Xx\s \in \D(A)$ and
$AXx\s  +XA\s x\s= -Q x\s.$
\end{definition}

We recall the following observation 
from \cite{Vu}; since our setting is slightly different
we include a proof.

\begin{proposition}\label{prop:identity-for-X}
If $X$ is a positive symmetric solution of the equation
\eqref{eq:lyapunov}, then for all $t>0$ we have
\begin{equation}
\label{eq:identity-for-X}
X - S(t) X  S\s(t)  = Q_t.
\end{equation}
\end{proposition}
\proof
From \eqref{eq:lyapunov}
we have, for $x\s,y\s\in \D(A\s)$,
$$
\begin{aligned}
 \lb Q_t x\s, y\s\rb 
 & = -\int^t_0 \lb S(s)(AX+XA\s)S\s(s) x\s, y\s\rb \,ds
\\ & = -\int^t_0 \frac{d}{ds} \lb S(s)XS\s(s)x\s, y\s\rb\,ds
= \lb Xx\s, y\s\rb - \lb S(t)XS\s(t)x\s, y\s\rb.
\end{aligned}
$$
Since $\D(A\s)$ is weak$\s$-dense in $E\s$, it follows that
$ Q_t x\s = Xx\s - S(t) X  S\s(t)x\s$ for all $x\s\in\D(A\s)$.
Finally, since both $Q_t$ and $X$ are positive and symmetric,
and therefore weak$\s$-to-weakly continuous, 
it follows from this that
$Q_t x\s = Xx\s - S(t) X  S\s(t)x\s$ for all $x\s\in E\s$.
\qed

After these preparations we can state and prove our main result about
the Liapunov equation.
Under somewhat more restrictive conditions,
this result was proved in \cite{Za} for the case when $E$ is a Hilbert space; see also \cite[Theorem 11.7]{DZ}. 

\begin{theorem} \label{thm:lyapunov}
The following assertions are equivalent:

\begin{enumerate}
\item Equation \eqref{eq:lyapunov} has a positive symmetric solution;
\item Hypothesis $(\HQinfty)$ holds.
\end{enumerate}
If these equivalent conditions are satisfied, the operator $Q_\infty$ is
a positive symmetric solution of \eqref{eq:lyapunov}, which is minimal
in the sense that 
if $R$ is another positive symmetric solution of \eqref{eq:lyapunov}, 
then for all $x\s\in E\s$ we have
$$\lb \Q x\s,x\s\rb \le \lb Rx\s,x\s\rb.$$
\end{theorem}

\proof
(1)$\,\Rightarrow\,$(2): Suppose $R$ is a positive symmetric solution
of \eqref{eq:lyapunov}.
By \eqref{eq:identity-for-X}, for all $x\s\in E\s$ we have
$$
 \lb Q_t x\s, x\s\rb
 = \lb  R  x\s,x\s\rb - \lb  R  S\s(t)x\s, S\s(t)x\s\rb.
$$ 
Since by assumption $R$ is positive, this implies
$$
 \lb Q_t x\s, x\s\rb\le \lb R x\s,x\s\rb
$$
for all $x\s\in E\s$. Hence, $t\mapsto \lb Q_t x\s,x\s\rb$ is a 
bounded function. Since this function is 
also non-decreasing, it follows that the limit
$\lim_{t\to\infty} \lb Q_t x\s, x\s\rb$ exists for all $x\s\in E\s$
and we have
$$ \lim_{t\to\infty} \lb Q_t x\s, x\s\rb\le \lb Rx\s,x\s\rb.$$ 
By polarization,  the limit $\lim_{t\to\infty} \lb Q_t x\s, y\s\rb$ 
exists for all  $x\s, y\s\in E\s$.
We now define a linear operator $Q_\infty \in \L(E\s, E\s\s)$ by 
$$\lb y\s,  Q_\infty x\s\rb :=  \lim_{t\to\infty} \lb Q_t x\s, y\s\rb.$$
By the uniform boundedness theorem, $\Q$ is bounded. 
We claim that $\Q$ actually takes values in $E$. 
Indeed, for all $x\s\in E\s$ we have
$$
  \lb x\s, \Q x\s\rb  
 = \lim_{t\to\infty} \lb Q_t x\s, x\s\rb \le \lb R x\s,x\s\rb, \qquad
x\s\in E\s,
$$
and the claim follows from Proposition \ref{E-valued}.

\smallskip\noindent
(2)$\,\Rightarrow\,$(1): This is the content of Proposition
\ref{prop:lyapun}.
\qed

%

The semigroup $\SS$ is said to be $Q$-{\em symmetric} if
for all $t\ge 0$ we have $$S(t)Q = QS\s(t).$$
It is easy to check that the following assertions
are equivalent:

\begin{enumerate}
\item $\SS$ is $Q$-symmetric;
\item For all $x\s\in \D(A\s)$ we have $Qx\s\in \D(A)$ and
$AQx\s = QA\s x\s$.
\end{enumerate}
If Hypothesis ($\HQinfty$) holds, then these assertions are equivalent to:

\begin{enumerate}
\item[(3)] For all $t\ge 0$ we have $S(t)\Q = \Q S\s(t)$;

\item[(4)] For all $x\s\in \D(A\s)$ we have $\Q x\s\in \D(A)$ and
$A\Q x\s = \Q A\s x\s$.
\end{enumerate}
It follows from (4) and Liapunov equation that 
$$A\Q = \Q A\s = -\tfrac12 Q.$$

\begin{theorem} [$\HQinfty$] \label{thm:S_H(t)}
If $\SS$ is $Q$-symmetric, 
then $\SS$ restrict to a
selfadjoint and strongly stable $C_0$-semigroup of contractions on $H$.
\end{theorem}
\proof
For all $x\s\in \D(A\s)$
we have
$$ - \lb Qx\s, A\s x\s\rb = 2\lb A \Q x\s, A\s x\s\rb = 
2 \lb Q_\infty A\s x\s, A\s x\s\rb \ge 0 $$
since $Q_\infty$ is a positive operator.
Therefore by Theorem \ref{thm:contraction}, $\SS$ maps $H$ into itself and
the restricted semigroup $\SS_H$ is a $C_0$-semigroup of contractions on $H$. 

Selfadjointness of $\SS_H$ follows from
$$
\begin{aligned}
{[i\s x\s, S_H\s(t)i\s y\s]_H}
& = [S_H(t) i\s x\s, i\s y\s]_H
  = [i\s S\s(t) x\s, i\s y\s]_H
\\ & = \lb Qy\s, S\s(t) x\s\rb
 = \lb S(t)Qy\s, x\s\rb
\\ & 
 = \lb i S_H(t) i\s y\s, x\s\rb = [i\s x\s, S_H(t) i\s y\s]_H.
\end{aligned}
$$
It remains to prove strong stability of $\SS_H$. 
Fix $x\s\in E\s$. From
$$
\int^\infty_0 \n S_H\s(t)i\s x\s\n_H^2\,dt = \lb Q_\infty x\s,x\s\rb <
\infty
$$ 
and a standard argument
it follows that $$\lim_{t\to\infty} \n S_H(t)i\s x\s\n_H
= \lim_{t\to\infty} \n S_H\s (t)i\s x\s\n_H= 0.$$
By a density argument, this gives the strong stability
of $\SS_H$.
\qed

In the following result, which extends a result from \cite{CG3}, 
we do not assume Hypothesis $(\HQinfty$): 

\begin{proposition}\label{cor:selfadoint-C_0}
The following assertions are equivalent:

\begin{enumerate}
\item $\SS$ is $Q$-symmetric;
\item $H$ is $\SS$-invariant and $\SS_H$ is a selfadjoint semigroup on $H$;
\item $H$ is $\SS$-invariant and $\SS_H$ is a selfadjoint $C_0$-semigroup on $H$.
\end{enumerate}
\end{proposition}
\proof
(1)$\,\Rightarrow\,$(3): By rescaling $\SS$ 
we may  assume that $\SS$ is uniformly exponentially stable.
Then Hypothesis ($\HQinfty$) holds, and  
the assertion follows from Theorem \ref{thm:S_H(t)}.

(3)$\,\Rightarrow\,$(2): Trivial.

(2)$\,\Rightarrow\,$(1): For all $x\s,y\s\in E\s$,
$$
\begin{aligned}
\lb S(t) Qx\s, y\s\rb & = \lb i \,S_H(t)\, i\s x\s, y\s\rb
 = \lb i\, S_H\s(t)\,i\s x\s, y\s\rb
\\ & = \lb i\, i\s S\s(t) x\s, y\s\rb
 = \lb QS\s(t) x\s, y\s\rb.
\end{aligned}
$$
\qed


\section{Spectral gap conditions}
\label{sec:gap}


In this section we shall prove some results for 
the semigroup $\SS_\infty$, which will be applied in Section \ref{sec:OU2}
to obtain a necessary and sufficient condition for the
existence of a spectral gap for the
generator of the Ornstein-Uhlenbeck semigroup associated with $\SS$ and $Q$. 

\medskip
We start with a simple but useful lemma.

\begin{lemma}[$\HQinfty$]\label{lem:core}
The set $i_\infty\s(\D(A\s))$ is a core for $A_\infty\s$.
\end{lemma}
\proof
Since $\D(A\s)$ is weak$\s$-dense in $E\s$ and $i_\infty\s$ is weak$\s$-to-weakly continuous,
the set $i_\infty\s(\D(A\s))$ is weakly dense,
and hence dense, in $H_\infty$.

From $S_\infty\s(t) i_\infty\s = i_\infty\s \circ S\s(t)$ it follows that $i_\infty\s(\D(A\s))$
is $\SS_\infty\s$-invariant, 
and by another standard result from from the theory of $C_0$-semigroups 
\cite[Proposition II.1.7]{EN}
this implies that $i_\infty\s(\D(A\s))$ is a core for $A_\infty\s$.
\qed

The next result extends a result from \cite{CG2} to the Banach space setting.
 
\begin{lemma} [$\HQinfty$]
\label{thm:gap} 
Let $M>0$ be given. 
The following statements are equivalent: 

\begin{enumerate}
\item $ \lb Q_\infty x\s, x\s\rb \le M \lb Qx\s,x\s\rb$ for all $x\s\in E\s;$

\smallskip
\item $\n S_\infty(t)\n_{H_\infty} 
\le \exp\bigl(- \frac{t}{2M}  \bigr)$ for all $t\ge 0$;

\end{enumerate}
\end{lemma}

\proof
Before we start the proof we note that for all $x\s\in \D(A\s)$ we have
the identity
\begin{equation}\label{eq:identity}
\begin{aligned}
\lb Qx\s, x\s\rb 
& = -2 \lb x\s, A\Q x\s\rb.  
\\ & = -2 [i_\infty\s A\s x\s, i_\infty\s x\s]_{\H}
= -2 [A_\infty\s i_\infty\s x\s, i_\infty\s x\s]_{\H}.
\end{aligned}
\end{equation}

(1)$\,\Rightarrow\,$(2):
By \eqref{eq:identity},
 \begin{equation}\label{eq:identit-2}
 \n i_\infty\s x\s\n_{H_\infty}^2
= \lb Q_\infty x\s,x\s\rb
\le M  \lb Qx\s, x\s\rb 
= -2M [A_\infty\s i_\infty\s x\s, i_\infty\s x\s]_{\H}.
\end{equation}
Hence by Lemma \ref{lem:core}, for all $h\in \D(A_\infty\s)$ we obtain
\begin{equation}\label{eq:contr-sg}
[A_\infty\s h,h]_\H\le -\tfrac{1}{2M} \n h\n_\H^2.
\end{equation}
By standard results on contraction semigroups in Hilbert spaces, 
this is equivalent to (2).

(2)$\,\Rightarrow\,$(1): If (2) holds, then \eqref{eq:contr-sg} holds
for all $h\in \D(A_\infty\s)$. Taking $h = i_\infty x\s$ with $x\s\in \D(A\s)$,
it follows that \eqref{eq:identit-2} holds
for all $x\s\in  \D(A\s)$. It remains to prove that  \eqref{eq:identit-2}
holds
for all $x\s\in E\s$.

Let $x\s\in E\s $ be arbitrary and fixed, 
and let $(x_n\s)$ be a sequence in
$\D(A\s)$ converging to $x\s$ weak$\s$ in $E\s$. 
Then $i\s x_n\s \to i\s x\s$ weakly in $H$ and 
$i_\infty\s x_n\s \to i_\infty\s x\s$ weakly in $H_\infty$.
Choose a sequence $(y_n\s)$ consisting of convex
combinations of elements from $(x_n\s)$ for which
$i\s y_n\s\to i\s  x\s$ strongly in $H$
and $i_\infty
y_n\s\to i_\infty\s x\s$ strongly in $H_\infty$.

Applying \eqref{eq:identit-2} to $y_n\s$ 
and passing to the limit for $n\to\infty$ gives
$$\lb Q_\infty x\s,x\s\rb
= \n i_\infty\s x\s\n_{H_\infty}^2 
\le M \n i\s x\s\n_H^2
= M \lb Qx\s,x\s\rb
$$
for all $x\s\in E\s$, and we obtain (1).
\qed

Combined with Proposition \ref{prop:Hinfty-2} this implies:

\begin{corollary}[$\HQinfty$]\label{cor:H_t}
 If $H_\infty\subseteq H$, then for all $t>0$ we have $H_t = H_\infty$ 
with equivalent norms.
\end{corollary}

In case $H$ is invariant, we can say more. For the $Q$-symmetric case in
Hilbert spaces, the following result was obtained in \cite{CG3}. Let $A_H$ denote the generator of the semigroup $\SS_H$. 

\begin{theorem}[$\HQinfty$]\label{thm:gap-2}
If $\SS$ restricts to a $C_0$-semigroup on $H$, 
the following assertions are equivalent:
\begin{enumerate}
\item $\SS_H$ is uniformly exponentially stable;
\item $\SS_\infty$ is uniformly exponentially stable;
\item For some $t>0$ we have $H_t = H_\infty$ with equivalent norms;
\item For all $t>0$ we have $H_t = H_\infty$ with equivalent norms;
\item $H_\infty \subseteq H$.
\end{enumerate}
In this situation, the inclusion $H_\infty\subseteq H$ is dense.
\end{theorem}
\proof
(2)$\,\Leftrightarrow\,$(3): Recall that 
$\SS_\infty$ is uniformly exponentially stable if and only if
there exists $t_0>0$ such that
$\n S_\infty(t_0)\n_{\H} <1$. 
Therefore the equivalence follows from 
Proposition \ref{prop:Hinfty-2}.

(3)$\,\Rightarrow\,$(4): This follows from Theorem \ref{thm:H-invar}.

(4)$\,\Rightarrow\,$(5): Fix $t_0>0$. Then 
by Proposition \ref{prop:str-cont}, 
$\H = H_{t_0}\subseteq H$.

(5)$\,\Rightarrow\,$(2): 
This follows from Lemma \ref{thm:gap}.

(5)$\,\Rightarrow\,$(1) : 
By assumptions there is a constant $K$ such that
for all $x\s\in E\s$ we have
$$ \int^\infty_0 \n S_H(t) i\s x\s\n_H^2\,dt 
= \lb Q_\infty x\s,x\s\rb 
\le K \lb Qx\s,x\s\rb = K\n i\s x\s\n_H^2.
$$
Hence the map $i\s x\s \mapsto S_H(\cdot)i\s x\s$ 
has a bounded extension, of norm $\le \sqrt{K}$,
to a bounded operator from $H$ into $L^2([0,\infty);H)$.
But then the Datko-Pazy theorem \cite{Pa} implies (1). 

(1)$\,\Rightarrow\,$(5): 
For all $x\s\in E\s$ we may estimate
$$\begin{aligned}
\lb Q_\infty x\s,x\s\rb 
 & = \int^\infty_0 \n S_H(t) i\s x\s\n_H^2\,dt 
\\ & \le \int^\infty_0  Me^{-at} \,dt \cdot \n i\s x\s\n_H^2
  = Ma^{-1} \lb Qx\s,x\s\rb.
\end{aligned}
$$
This gives (5).

The final assertion follows from Proposition 
\ref{prop:H-cap-Hinfty}. Alternatively, one could observe that
by Proposition 
\ref{prop:str-cont}, the inclusions $H_t\subseteq H$ are dense.
Therefore the result follows from the fact that $H_t=\H$ with equivalent
norms.
\qed

In view of Lemma \ref{thm:gap},
it seems natural to ask whether assertions (2) and (5) in Theorem \ref{thm:gap-2}
are always equivalent (i.e., even when $H$ fails to be invariant). 
The following example shows that this is not the case.

\begin{example}\label{ex:Hinfty}
Let $E = \R^2$, let 
$$\displaystyle
Q \ =\  \left(
\begin{array}{cc}
0 & 0 \\
0 & 1 
\end{array}
\right)
$$ 
and let the semigroup $\SS$ on $E$ be given by
$$\displaystyle 
S(t) \ = \ e^{-t}\left(
\begin{array}{cc}
1 & t \\
0 & 1 
\end{array}
\right)
.$$ 
An easy computation gives
$$
\Q \ =\  \frac14 \left(
\begin{array}{cc}
1 & 1\\
1 & 2 
\end{array}
\right),
$$
cf. \cite[Example 4.3]{fu}.
This matrix is invertible, so $H_\infty = E = \R^2$ (with equivalent norms).
On the other hand, $H$ is the one-dimensional subspace of $E$ spanned 
by the vector 
$\displaystyle\left( \begin{array}{c} 0\\ 1 \end{array}\right)$.
It follows that $$H_\infty\not\subseteq H.$$
On the other hand it is clear that $\SS_\infty$ is strongly stable,
hence uniformly exponentially stable since $H_\infty$ 
is finite-dimensional.
 
There is no contradiction with Theorem \ref{thm:gap}: the point is that 
there exists no $\omega>0$
such that
$\n S_\infty(t) \n_{H_\infty} \le  e^{-\omega t}$ for all $t\ge 0$.
This can be checked by the following direct computation that
will be useful in the next section as well. 

We use the simple fact that 
$\n S_\infty(t)\n^2$ is the largest eigenvalue of $S_\infty(t)S_\infty\s(t)$.
Noting that $i_\infty\s$ is a surjection from $E$ onto $\H$, 
the number $\lambda(t)$ is an eigenvalue of $S_\infty(t)S_\infty\s(t)$
if and only if there exists a vector $x\s\in E$ such that for all $y\in E\s$ we
have
$$ [S_\infty(t)S_\infty\s(t)\i\s x\s, \i\s y\s]_{\H} = 
\lambda(t)[\i\s x\s, \i\s y\s]_{\H},
$$
or equivalently,
$$ \lb S(t)\Q S\s(t)\s x\s, y\s\rb = \lambda(t)\lb \Q x\s, y\s\rb.$$
Thus we have to solve the equation 
$$ \det \, \bigl[S(t)\Q S\s(t) - \lambda(t) \Q\bigr] \ =\ 0.$$
An elementary computation gives
$$ S(t)\Q S\s(t) - \lambda(t) \Q
= \frac14\left(
\begin{array}{cc}
e^{-2t}\bigl(2t^2+2t+1  \bigr) -\lambda(t) &e^{-2t}\bigl(2t+1  \bigr) -\lambda(t)   \\
e^{-2t}\bigl(2t+1  \bigr) -\lambda(t) &2e^{-2t} -2\lambda(t)   \\
\end{array}\right)
$$
from which we deduce that
$$ \lambda_{\pm}(t)  = e^{-2t} \left(t \pm \sqrt{t^2+1}\right)^2.$$
We finally obtain
$$\n S_\infty(t)\n =e^{-t} \bigl(t + \sqrt{t^2+1}\bigr).$$
Denoting the right hand side by $f(t)$, we have $f(0)=1$,
$\lim_{t\to\infty}f(t) = 0$ mono\-tonously, and
$f'(0)= 0$. Clearly, a function $f(t)$ with these properties
cannot be dominated by a negative exponential $e^{-\omega t}$.
\end{example}

\section{The Ornstein-Uhlenbeck semigroup in $C_b(E)$}
\label{sec:OU1}

Positive symmetric operators from $E\s$ into $E$ arise naturally as
the covariance operators of  Gaussian Borel measures on $E$.
However, not every positive symmetric operator is a Gaussian
covariance operator,
and for this reason we will frequently consider the following
hypothesis:
\begin{itemize}

\item ($\Hmu$): $E$ is separable,
and for all $t>0$ the operator $Q_t$ is the covariance of a centred
Gaussian Borel measure $\mu_t$ on $E$. \smallskip
\end{itemize}

Whenever it is convenient we further put $\mu_0:=\delta_0$, the Dirac measure
concentrated at the origin.

The separability in Hypothesis ($\Hmu$)
is added in order avoid certain measure theoretic
complications.

If $E$ is a separable real Hilbert space, then $Q_t$ is a Gaussian covariance
if and only $Q_t$ is a trace class operator, and this happens if and only if the
inclusion $i_t:H_t\embed E$ is a Hilbert-Schmidt operator.

The relevance of Hypothesis ($\Hmu$) is explained by the following result
from \cite{BN}:

\begin{proposition}[$\Hmu$]\label{prop:sACP}
Let $A$ denote the generator of the semigroup $\SS$.
The stochastic evolution equation
\begin{equation}\label{sACP}
\begin{aligned}
dX(t) & = AX(t)\,dt + dW_H(t), \qquad t\ge 0,\\
 X(0) & = x,
\end{aligned}
\end{equation}
has a unique weak solution $\{X(t,x)\}_{t\ge 0}$
if and only Hypothesis $(\Hmu)$ holds. In this situation
the process $\{X(t,x)\}_{t\ge 0}$ is Gaussian.
For all $t > 0$ we have $X(t,x) = S(t)x + X(t,0)$ almost surely, and
the distribution of $X(t,0)$ equals $\mu_t$.
\end{proposition}

Assuming Hypothesis ($\Hmu$), we define
the transition semigroup $\PP = \{P(t)\}_{t\ge 0}$
of $\{X(t,\cdot)\}_{t\ge 0}$ on the space $B_b(E)$
of all bounded Borel functions on $E$  by
$$ P(t)f(x) = \E \bigl(f(X(t,x)) \bigr)
= \int_E f(S(t)x+y)\,d\mu_t(y), \qquad t\ge 0, \ x\in E.$$
The semigroup $\PP$ is contractive on $B_b(E)$
and it maps $C_b(E)$,
the space of all bounded
continuous functions on $E$,
into itself.

In general, the semigroup $\PP$ is not strongly continuous on $C_b(E)$,
and not even on its closed subspace $BUC(E)$ of all bounded and uniformly
continuous functions on $E$.

We will show next that $\PP$ is strongly continuous
on $C_b(E)$ endowed with the {\em mixed topology} which is defined as the finest locally
convex topology on $C_b(E)$ that agrees on every norm-bounded set 
with the topology of uniform convergence on compact sets.    For Hilbert spaces $E$,
this fact was proved in \cite{GK}.
By the results in \cite{Se}, this definition agrees with the one in
\cite{GK}. Clearly,
$$\tau_{\rm uniform-on-compacts}\ \subset\ \tau_{\rm mixed}\ \subset\ \tau_{\rm uniform}.$$
We have the following characterization of sequential convergence in the mixed
topology: 
{\em a sequence $(f_n)$ in $C_b(E)$ converges to $f\in C_b(E)$ if and only if}
\begin{enumerate}
\item $\sup_n \n f_n\n_\infty<\infty$;
\item $\lim_{n\to \infty} f_n = f$ {\em uniformly on compact subsets of} $E$.
\end{enumerate}
We will also need the fact  that the dual space
$\left(C_b(E),  \tau_{\rm mixed}\right)^*$  can be identified in the natural way
with the space of finite Borel measures on $E$ \cite{fremlin}.

For more information about the mixed 
topology we refer the interested reader to 
the papers \cite{Se, Wh, Wi} and the references therein.

\begin{theorem}[$\Hmu$]\label{thm:mixed}
The semigroup $\PP$ is strongly continuous on $C_b(E)$ in its mixed topology.
\end{theorem}
\begin{proof}
Following the arguments of \cite{GK}, we see that it suffices to prove that
for all $f\in C_b(E)$ and all compact subsets $K\subseteq E$ we have 
$$\lim_{t\downarrow 0}\,\left( \sup_{x\in K} |P(t)f(x) -f(x)|\right) = 0.$$
For Hilbert spaces $E$, this can be proved easily by probabilistic arguments.
Here we give a direct, analytical proof.

Fix $f\in C_b(E)$ and $K\subseteq E$ compact. We may assume that $K$ is convex. 
As was observed in \cite{Ne}, we have
weak convergence $\mu_t\to\mu_0 = \delta_0$, the Dirac measure
concentrated at $0$. Fixing an arbitrary $\e>0$, 
by tightness we may choose
a compact set $L$ in $E$ such that $\mu_t(L)\ge 1-\e$ for all
$t\in [0,1]$. We may assume that $L$ is convex. 
Keeping in mind that $\mu_0 = \delta_0$, 
we  necessarily have $0\in L$.

For all $t\ge 0$ and $x\in E$ we have
\begin{equation}\label{split}
\begin{aligned}
\ &| P(t)f(x) - f(x)| 
\\ & \qquad \le \int_E |f(S(t)x+y) - f(x+y)|\,d\mu_t(y)
+ \int_E |f(x+y) - f(x)|\,d\mu_t(y).
\end{aligned}
\end{equation}
We will estimate the two integrals on the right hand side separately.
For all $x\in K$ and $t\in [0,1]$ we have
$$
\int_E |f(S(t)x+y) - f(x+y)|\,d\mu_t(y)
\le 2\e \,\n f\n_\infty + \int_L |f(S(t)x+y) - f(x+y)|\,d\mu_t(y).
$$
By the strong continuity of $\SS$, which is uniform on compact sets,
and the uniform continuity of $f$ on the
compact set $\{S(t)x+y: \ (t,x,y)\in [0,1]\times K\times L\}$ 
we may choose $0<t_0\le 1$ so small that
$$\sup_{x\in K, \, y\in L} |f(S(t)x+y)-f(x+y)| <\e, \qquad t\in [0,t_0]. $$
Thus, for all $t\in [0,t_0]$ we obtain
\begin{equation}\label{est1}
\sup_{x\in K}\int_E |f(S(t)x+y) - f(x+y)|\,d\mu_t(y)
\le 2\e\, \n f\n_\infty + \e.
\end{equation}
Next we estimate the second integral on the right hand side of \eqref{split}.
As above, for all $x\in K$ and $t\in [0,1]$ we have
\begin{equation}\label{est2a} 
\int_E |f(x+y) - f(x)|\,d\mu_t(y) \le 2\e\,\n f\n_\infty + 
\int_L |f(x+y) - f(x)|\,d\mu_t(y).
\end{equation}
Hence it remains to show that
\begin{equation}\label{remains}
\lim_{t\downarrow 0}\left(\sup_{x\in K}  \int_L |f(x+y) -
f(x)|\,d\mu_t(y)\right) = 0.
\end{equation}
The restriction of $f$ to $K+L$ being uniformly continuous, we 
introduce its modulus of continuity,
$$
\eta(\delta) := \sup\{ |f(z)-f(z')|: \ z,z'\in K+L, \
\n z-z'\n\le \delta\}.
$$
Then, recalling that $0\in L$,
$$ \sup_{x\in K} \int_L |f(x+y) - f(x)|\,d\mu_t(y)\le
\int_L \eta(\n y\n)\,d\mu_t(y).
$$
The function $\zeta(y) :=  \eta(\n y\n)$ is bounded, nonnegative, and
continuous on $L$. By the Tietze-Urysohn extension
theorem \cite[Theorem 2.1.8]{Eng}, it can be extended to a bounded,
nonnegative, and continuous function $\zeta$ on all of $E$.
The weak convergence $\mu_t\to\delta_0$ then implies
$$\limsup_{t\downarrow 0}  \int_L \eta(\n y\n)\,d\mu_t(y)
\le \lim_{t\downarrow 0}  \int_E \zeta(y)\,d\mu_t(y)
= \zeta(0) = 0.
$$
This proves \eqref{remains}.
\end{proof}

In the remainder of this section we will always consider $\PP$ as a
strongly continuous semigroup on $C_b(E)$ in its mixed topology.
The {\em infinitesimal generator} $(L,\D(L))$ of $\PP$ is defined by 
$$
\begin{aligned}
\D(L) & = \left\{ f\in C_b(E): \ \lim_{t\downarrow 0} \frac{P(t)f - f}{t} \
\hbox{ exists }\right\}, \\
Lf & = \lim_{t\downarrow 0} \frac{P(t)f - f}{t} \quad (f\in \D(L)),
\end{aligned}
$$
where the limits are taken with respect to the mixed topology. In a
similar way we define the weak generator: we say that
$\phi\in \D\left(L_w\right)$ if there exists a (necessarily unique)
function $f\in C_b(E)$ such that
\begin{equation}\lim_{t\downarrow 0}\int_E\frac {P(t)\phi (x)-\phi (x)}{t}\,d\nu(x)
= \int_E f(x)\,d\nu(x)\label{a2}\end{equation}
for each  finite Borel measure $\nu$ on $E$; then we
define $L_w\phi :=f$.

\begin{proposition}[$\Hmu$]\label{ts2}
We have $L=L_w$.
\end{proposition}
\begin{proof}
The proof is a simple modification of the proof of
\cite[Corollary 1.2]{Pa}.
Obviously $L\subseteq L_w$. Let, conversely, $\phi\in\D\left(L_w\right
)$. Then,  as in
\cite{Pa}, for each finite Borel measure $\nu$ on $E$ we find
\[\int_E P(t)\phi (x)-\phi (x)\, d\nu(x)=\int_E\left(\int_
0^t P(s)L_w\phi (x)\,ds\right)\,d\nu(x).\]
Since this holds for all finite Borel measures $\nu$ and the functions $x\mapsto \phi(x)$, $x\mapsto P(t)\phi(x)$ and $x\mapsto 
\int^t_0 P(s)L_w\phi(x)\,ds$
are continuous we obtain, for all $x\in E$,
\[P(t)\phi(x) -\phi(x) =\int_0^t P(s)L_w\phi(x)\, ds.\]
Hence if $K\subseteq E$ is compact, then
$$
\begin{aligned}
 & \lim_{t\downarrow 0}\left(
\sup_{x\in K} \left[\frac {P(t)\phi(x) -\phi(x)}{t} - L_w\phi(x)\right]\right)
\\ & \qquad \qquad =
\lim_{t\downarrow 0}\left(
\sup_{x\in K}\left[ \frac1t \int_0^t P(s)L_w\phi(x) - L_w\phi(x)\, ds\right]\right)
=0
\end{aligned}
$$
where in the last step we used that $\lim_{t\downarrow 0} P(t)L_w\phi(x) = L_w\phi(x)$
uniformly on $K$. By the definition of $L$ and the afore-mentioned criterium for sequential
convergence in the mixed topolgy, this shows that $\phi\in\D(L)$.
This  concludes the proof.
\end{proof}

As an immediate corollary we have the following result,
which shows that our definition of $(L,\D(L))$ agrees with 
the one in \cite{priola}.

\begin{corollary}\label{ts3}
We have $\phi\in \D(L)$ if and only if the following two conditions hold:

\begin{enumerate}
\item $\displaystyle\limsup_{t\downarrow 0}
\frac{\left\|P(t)\phi -\phi\right \|}{t}<\infty$; 
\item There exists a function $f\in C_b(E) $ such that for all $x\in E$,
$$\lim_{t\downarrow 0}\frac {P(t)\phi (x)-\phi (x)}{t}=f(x).$$
\end{enumerate}
In this situation we have $L\phi = f$.
\end{corollary}

Our next aim is to obtain an explicit representation of $L$ on a suitable core.

\begin{lemma}[$\Hmu$]\label{prop1}
Let $\{X(t,x_0)\}_{t\ge 0}$ be the unique weak solution of \eqref{sACP}.
Let $\phi\in C^2(\R^d)$,  $x_1\s, \hdots, x_d\s\in \D(A\s)$, and $x_0\in E$. 
Then for all $t\ge 0$ the following identity holds almost surely:
$$
\begin{aligned}
\phi(\lb & X(t,x_0),  x_1^*\rb,\hdots, \lb X(t,x_0), x_d^*\rb)  
\\ & = \phi(\lb x_0, x_1^*\rb,\hdots, \lb x_0, x_d^*\rb) 
\\ &  \quad + \sum_{j=1}^d\int^t_0 
    \frac{\partial \phi}{\partial x_j}(\lb X(s,x_0),x_1^*\rb,\hdots,
    \lb X(s,x_0), x_d^*\rb)
    \, d W_H(s) i^* x_j^*
\\ & \quad   + \sum_{j=1}^d \int^t_0 
      \frac{\partial \phi}{\partial x_j}(\lb X(s,x_0),x_1^*\rb,\hdots, \lb X(s,x_0), x_d^*\rb)
      \lb X(s,x_0),A\s x_j\s\rb \,ds
\\ &\quad   + \frac12\sum_{j=1}^d \sum_{k=1}^d [i\s x_j\s,i\s x_k\s]_H
\int^t_0
\frac{\partial^2 \phi}{\partial x_j\partial x_k}(\lb X(s,x_0),x_1^*\rb,
\hdots, \lb X(s,x_0), x_d^*\rb) \,ds.
\end{aligned}$$
\end{lemma}
\proof
For $j=1,\dots, d$ and $t\ge 0$ we define
$$
 \xi_j(t) := \lb x_0,x_j\s\rb + \int^t_0 \lb X(s,x_0), A\s x_j\s\rb\,ds 
 + W_H(t)i\s x_j\s.
$$
By the definition of a weak solution, for all $t\ge 0$ we have
$\xi_j(t) = \lb X(t,x_0), x_j\s\rb$ almost surely, 
so $\{\xi_j(t)\}_{t\ge 0}$ is a modification of
the process $\{\lb X(t,x_0), x_j^*\rb\}_{t\ge 0}$.
Since almost surely, the trajectories of $\{\lb X(t,x_0),A\s x_j\s\rb\}_{t\ge 0}$
are locally integrable, we see that almost surely the trajectories
of the process
$ \{V_j(t)\}_{t\ge 0}$ defined by
\begin{equation}\label{V_j}
V_j(t): = \int^t_0 \lb X(s,x_0), A\s x_j\s\rb\,ds.
\end{equation}
are continuous and locally of bounded
variation.
By redefining the random variables $V_j(t)$ to be $0$
on a common null set we obtain a modification
of $\{\lb X(t,x_0), x_j^*\rb\}_{t\ge 0}$, still denoted
by  $\{\xi_j(t)\}_{t\ge 0}$,
From the representation
$$ \xi_j(t) = \lb x_0, x_j\s\rb + V_j(t) + M_j(t)$$
where $M_j(t) := W_H(t)i\s x_j\s$, we see
that $\{\xi_j(t)\}_{t\ge 0}$ is a
continuous semimartingale.  
Define $F:\R^d\times\R^d\to\R$ by 
$$ F(u,v) := \phi\bigl((\lb x_0,x_1\s\rb,\hdots,\lb x_0,x_d\s\rb)+u+v\bigr).$$
Put
$$
\xi(t)  := (\xi_1(t),\dots,\xi_d(t)),\
  M(t)  :=(M_1(t),\hdots,M_d(t)), \
  V(t)  :=(V_1(t),\hdots,V_d(t)).  
$$
By the It\^o formula \cite[Theorem 5.10]{CW} almost surely we have, for all 
$t\ge 0$,
$$
\begin{aligned}
\phi(\xi(t))& -\phi(\xi(0))
  = F(M(t),  V(t)) - F(M(0), V(0)) 
\\ & = \sum_{j=1}^d \int^t_0 
     \frac{\partial F}{\partial u_j}(M(s), V(s))\,dM_j(s)
\\ & \qquad + \sum_{j=1}^d\int^t_0 
     \frac{\partial F}{\partial v_j}(M(s), V(s))\,dV_j(s)
\\ & \qquad +
\frac12\sum_{j=1}^d\sum_{k=1}^d[i\s x_j\s,i\s x_k\s]_H \int^t_0
\frac{\partial^2 F}{\partial u_j\partial v_k }(M(s), V(s))\,ds
\end{aligned}
$$
\begin{equation}\label{ito}
\begin{aligned}
 & = \sum_{j=1}^d\int^t_0 
    \frac{\partial \phi}{\partial x_j}(\xi(s)) 
    \, d W_H(s) i^* x_j^*
\\ & \qquad   + \sum_{j=1}^d \int^t_0 
      \frac{\partial \phi}{\partial x_j}(\xi(s))
      \lb X(s,x_0),A\s x_j\s\rb \,ds
\\ &\qquad   + \frac12\sum_{j=1}^d \sum_{k=1}^d [i\s x_j\s,i\s x_k\s]_H
\int^t_0
\frac{\partial^2 \phi}{\partial x_j\partial x_k}(\xi(s)) \,ds,
\end{aligned}
\end{equation}
where we used \eqref{V_j} and the fact  
that the mutual quadratic variation of $M_j(s)$ and $M_k(s)$
equals $s\cdot[i\s x_j\s, i\s x_k\s]_H$.

We claim that almost surely,
\begin{equation}\label{cont_vers}
\xi_j(s) = \lb X(s,x_0),  x_j\s\rb
\quad\hbox{for almost all } s\ge 0. 
\end{equation}
To see this, note that $\{\xi(t)\}_{t\ge 0}$ is  
progressively measurable,  being a process with continuous trajectories.
Also, $\{\lb X(t,x_0),  x\s\rb \}_{t\ge 0}$ 
is progressively measurable, being predictable.
The claim follows from Fubini's theorem.

The proposition now follows by combining \eqref{ito} and \eqref{cont_vers}.
\qed

We will now identify a suitable core for $L$ consisting of
cylindrical functions satisfying the assumptions of Lemma \ref{prop1}.
To this end let us define 
$\F = \FF$ as the space of all functions $f:E\to\R$ of the 
form
\begin{equation}
\label{eq:cyl}
f(x) = \phi(\lb x,x_1\s\rb, \hdots,\lb x,x_d\s\rb)
\end{equation}
for some $d\ge 1$, 
with $x_j\s\in\D(A\s)$ for
all $j=1,\hdots,d$ and $\phi\in C_b^2(\R^d)$.
Let 
$$\F_0 = \{f\in \F: \lb\, \cdot\,, A\s Df(\cdot)\rb\in C_b(E)\}.$$
For $f\in\F_0$ we define $L_0 f\in C_b(E)$ by 
\begin{equation}\label{L0}
L_0 f(x) := \tfrac12 \,\hbox{trace}\, D_H ^2 f(x) + \lb x, A\s Df(x)\rb,
\qquad x\in E.
\end{equation}
Here $Df:E\to E\s$ is the Fr\'echet derivative of $f$,
$$
Df(x) = \sum_{j=1}^d \frac{\partial\phi}{\partial x_j}(\lb
x,x_1\s\rb,\hdots,\lb x,x_d\s\rb)\otimes x_j\s
$$
and $D_Hf:E\to H$ is defined by 
\begin{equation}\label{defHQ}
D_{H}f(x)  = \sum_{j=1}^d \frac{\partial\phi}{\partial x_j}(\lb
x,x_1\s\rb,\hdots,\lb x,x_d\s\rb)\otimes i\s x_j\s.
\end{equation}
In a slightly different setting, the space $\F_0$ was introduced first by
Cerrai and Gozzi \cite{CG}; see also \cite{GK}.
Extending the results from these papers to the Banach space setting, 
we will show in a moment that $\F_0$
is a core for the generator $L$ and that 
$Lf=L_0 f$ for $f\in\F_0$.

\begin{theorem}[$\Hmu$]
$\F_0$ is a core for $L$, and for all $f\in\F_0$ we have $Lf = L_0 f$.
\end{theorem}
\begin{proof}
We will show first that $L_0\subseteq L$.
Clearly, $L_0f\in C_b(E)$ for any $f\in  F_0$.
Let $f(x) = \phi(\lb x,x_1\s\rb, \hdots,\lb x,x_d\s\rb)$ with $x_j\s\in\D(A\s)$ for
all $j=1,\hdots,d$ and $\phi\in C_b^2(\R^d)$.
First we note that
$$D_H ^2 \phi(\lb x,x_1\s\rb, \hdots,\lb x,x_d\s\rb) =
\sum_{j=1}^d \sum_{k=1}^d
\frac{\partial^2 \phi}{\partial x_j\partial x_{k}}
(\lb x,x_1\s\rb, \hdots,\lb x,x_d\s\rb) \otimes (i\s x_j\s\otimes
B\s x_{k}\s).$$
Now we apply Lemma \ref{prop1} and take on both sides the expectation.
This gives
$$
\begin{aligned}
\frac{1}{t} & \bigl( P(t)f(x) - f(x)\bigr)
\\ & =\E\,\frac{1}{t} \bigl(f(X(t,x)) - f(x) \bigr)
\\ & =\E\,  \frac1t \int^t_0
      \sum_{j=1}^d \frac{\partial \phi}{\partial x_j}(\lb X(s,x),x_1^*\rb,\hdots,
      \lb X(s,x), x_d^*\rb)\lb X(s,x),A\s x_j\s\rb \,ds
\\ & \qquad + \E\, \frac{1}{2t}
\sum_{j=1}^d\sum_{k=1}^d[i\s x_j\s,i\s x_{k}\s]_H \int^t_0
\frac{\partial^2\phi}{\partial x_j\partial x_{k}}
(\lb X(s,x),  x_1\s\rb,\hdots,\lb X(s,x),x_d\s\rb)\,ds
\\ & = \frac1t \int^t_0
       P(s)\bigl(\bigl\lb \,\cdot\,, A\s Df(\cdot)\bigr\rb\bigr)(x)  \,ds
       + \frac{1}{2t}\int^t_0\,P(s)\bigl(\hbox{trace}\,
D_H ^2 f\bigr)(x) \,ds.
\end{aligned}
$$
The assumption $f\in\F_0$ implies that the functions
$$x\mapsto \lb x, A\s Df(x)\rb \ \ \hbox{ and } \ \ x\mapsto
\tfrac{1}{2}\, \hbox{trace}\, D_H^2 f(x)$$
belong to $C_b(E)$. Therefore by Theorem \ref{thm:mixed}
we can pass to the limit $t\downarrow 0$ and by using Corollary \ref{ts3}
we obtain $L_0\subseteq L$ .
By \cite[Lemma 4.6]{GK} $\F_0$ is
dense in $C_b(E)$,  and since $\F_0$ is also $\PP-$invariant, $\F_0$ is
a core for $L$ by \cite[Lemma 4.7]{GK}; {\em mutatis mutandis},
the proofs of these results extend to the Banach space case.
\end{proof}

\subsection{The strong Feller property}
 We say that 
$\PP$ has the
{\em strong Feller property} if for every $t>0$, $P(t)$ maps
$B_b(E)$ into $C_b(E)$.
 We start with a characterization of the strong Feller property in terms of the mixed topology on $C_b(E)$.
Although a more general version of the following 
result below appears to be known to specialists,
it is not easily available and for the convenience of reader 
we include a straightforward 
proof for the case of
Ornstein-Uhlenbeck semigroups.

\begin{proposition}\label{sfpmixed}
The following conditions are equivalent:
\begin{enumerate}
\item The semigroup $\PP$ has the strong Feller property.
\item For each $t>0$, the mapping
\[ x\mapsto \mu_t(x,\cdot) \]
is continuous from $E$ into 
$\left(C_b(E), \tau_{{\rm mixed}}\right)^*$ with the variation norm.
\item For each $t>0$, the operator
\[P(t):  \left(C_b(E), \tau_{{\rm uniform}}\right)\to \left(C_b(E),
\tau_{{\rm mixed}}\right)\]
is compact.
\end{enumerate}
\end{proposition}
\begin{proof}
(1)$\,\Rightarrow\, $(2): Assume that (1) holds. 
As we have already mentioned, the dual space
$ \left(C_b(E), \tau_{{\rm mixed}}\right)^*$ may be identified as the space of
finite Borel
measures on $E$. 
By \cite[Theorem 9.19]{DZ}, for each $t>0$ there eixsts $c_t\ge 0$ such that
$$
 \left\|\mu_t(x,\cdot)-\mu_t(y,\cdot)\right\|_{{\rm var}}
\le c_t\n x-y\n,\quad x,y\in E,
$$
and (2) follows. 

(2)$\,\Rightarrow\, $(1):
If (2) holds, then by definition of the variation norm,
\[\left|P(t)\phi(x)-P(t)\phi(y)\right|\le 
\left\|\mu_t(x,\cdot )-\mu_t (y,\cdot )\right\|_{{\rm var}} \n \phi\n,\label{bg1}\]
which implies (1).
 
(2)$\,\Leftrightarrow\, $(3): Assume that (2) holds and let
\[B=\left\lbrace \phi\in C_b(E): \|\phi\|\le 1\right\rbrace.\]
The set $P(t)B$  is bounded and by (\ref{bg1}) uniformly
equicontinuous on compacts.
Therefore $P(t)B$ is relatively compact by \cite[Theorem 2.4]{GK} and 
(3) follows. 

(3)$\,\Leftrightarrow\,$(2): 
Assume now that (3) holds. Then
there exists a sequence $\left( \phi_n\right)$ in $B$ such that
\[\left\|\mu_t(x,\cdot )-\mu_t(y,\cdot )\right\|_{{\rm var}}=\lim_{n\to\infty}
\left|P(t)\phi_n(x)-P(t)\phi_n(y)\right|.\]
and since the set $P(t)B$
is relatively compact we may assume that $\limn P(t)\phi_n = \psi$ in the mixed topology. Hence, $\psi$ is continuous and (b) follows.
\end{proof}

It is well known that $\PP$ has the strong Feller property if and only if
the pair $(\SS,i)$ is null controllable \cite{DZ, Ne}.
Under the assumption that $H$ is $\SS$-invariant, from Theorem \ref{thm:null}
we thus obtain the following
explicit necessary and sufficient condition
for the strong Feller property. It
extends a previous result from \cite{CG3}.

\begin{theorem}[$\Hmu$] If $\SS$ restricts to a $C_0$-semigroup on $H$, then the
following assertions are equivalent:

\begin{enumerate}
\item $\PP$ has the strong Feller property;
\item $S(t)E \subseteq H$ for all $t>0$.
\end{enumerate}
\end{theorem}

\section{The Ornstein-Uhlenbeck semigroup in $L^2(E,\mu_\infty)$}
\label{sec:OU2}

In order to be able to 
study the semigroup $\PP$ in an $L^2$-context, we consider the
following hypothesis:

\begin{itemize}

\item ($\Hmuinfty$): $E$ is separable, 
Hypothesis ($\HQinfty$) holds, and the operator $Q_\infty$ 
is the covariance of a centred
Gaussian Borel measure $\mu_\infty$ on $E$.

\end{itemize}
By \eqref{eq:incl} and a standard tightness argument, ($\Hmuinfty$) implies ($\Hmu$).
If $E$ is a 
Hilbert space, then ($\Hmuinfty$) holds if and only if
($\Hmu$) holds and
$$\sup_{t>0}\
\bigl(\hbox{trace}\,Q_t\bigr) \, < \, \infty,$$
in particular, if $\SS$ is uniformly exponentially stable
\cite[Theorem 11.11]{DZ}. 
Extensions of these results to Banach spaces not containing a closed subspace
isomorphic to $c_0$ have been obtained in \cite{NW}.

In general it is not true that ($\Hmu$) and ($\HQinfty$) 
imply ($\Hmuinfty$), as is shown by the following example.

\begin{example}\label{diagonal}
Let $E=\ell^2$; we identify $E$ and its dual in a natural way.
The standard unit basis of $E$ will be denoted by $(x_n)_{n=1}^\infty$.

Let $(q_n)_{n=1}^\infty$ be a bounded sequence of strictly positive real
numbers and 
define $Q\in \L(E)$ by $$Qx_n := q_n x_n.$$
Let $(a_n)_{n=1}^\infty$ be a sequence of
strictly positive real numbers and define the operator $(A,\D(A))$ by
$$Ax_n := -a_n x_n$$
with maximal domain. 
Then $A$ generates a strongly stable $C_0$-semigroup $\SS$ on $E$, given by 
$S(t)x_n = e^{-a_nt}x_n.$

It is easy to check that
$$H = \left\{(b_n)_{n=1}^\infty\in E:\ 
 \sum_{n=1}^\infty \frac{1}{q_n} b_n^2 < \infty\right\}.
$$
For all $t>0$ we have 
$\displaystyle Q_t x_n = \frac{q_n}{2a_n} \bigl(1-e^{-2a_n t}\bigr)x_n,$
hence 
$$
 H_t 
= \left\{(b_n)_{n=1}^\infty\in E:\ 
\sum_{n=1}^\infty \frac{2a_n}{q_n}\bigl(1-e^{-2a_nt}\bigr)^{-1} b_n^2 < \infty\right\}.
$$
Hypothesis ($\HQinfty$) holds if and only if
$\displaystyle \sup_{n\ge 1} \ \frac{q_n}{a_n} \ < \, \infty.$
In this case we have
$\displaystyle \Q x_n = \frac{q_n}{2a_n}x_n$ and 
$$ \H = \left\{(b_n)_{n=1}^\infty\in E:\ 
 \sum_{n=1}^\infty \frac{a_n}{q_n} b_n^2 < \infty\right\}.
$$
Therefore, 
Hypothesis ($\Hmuinfty$) is satisfied if and
only if 
$\displaystyle
\sum_{n=1}^\infty \frac{q_n}{a_n} <\infty.
$

Let us now assume that the sequence $(a_n)_{n=1}^\infty$ is {\em bounded}.
Then there exists, for every $t>0$, a constant $M_t\ge 1$ such that
$$1\le  \frac{1-e^{-2a_nt}}{2a_n} \le M_t$$
for all $n\ge 1.$ From this it follows that
$$
 H_t
= \left\{(b_n)_{n=1}^\infty\in E:\ 
\sum_{n=1}^\infty \frac1{q_n} b_n^2 < \infty\right\}.
$$
Thus, for all $t>0$ we have $H_t=H$ up to an equivalent norm.
By computing the trace of $Q_t$ we see that  
 Hypothesis ($\Hmu$) is satisfied if and
only if $\displaystyle\sum_{n=1}^\infty {q_n} <\infty.$ 
For $q_n = 1/n^2$ and $a_n = 1/n$ we obtain an example
where ($\HQinfty$) and ($\Hmu$) hold, but not ($\Hmuinfty$). 
\end{example}
%

This example is interesting for another reason. It is shown in \cite{Ne3}
that $(\HQinfty$) implies that $A\s$ has no point spectrum in the closed
right half plane $\{z\in\C: \ \Re \,z\ge 0\}$, and that $A$
has no point spectrum on the imaginary axis if in addition
we assume that $\SS$ is uniformly bounded.
It may happen that $\sigma(A)\cap i\R$ is non-empty, however, even in the 
presence of Hypothesis ($\Hmuinfty$). For example, take
$q_n = 1/n^3$ and $a_n = 1/n$; then
$0\in\sigma(A)$
and Hypothesis ($\Hmuinfty$) holds.

Let us note that Example \ref{diagonal} can easily be extended to
$E = \ell^p$ ($1 \le p < \infty$) by using the fact \cite[Theorem V.5.6]{VTC} 
that a positive diagonal operator $(x_n)_{n=1}^\infty \mapsto (c_n x_n)_{n=1}^\infty$
 from $E\s$ to $E$ 
is a Gaussian covariance operator if and only if
$
\sum_{n=1}^\infty c_n^{p/2} <\infty.
$

\medskip
If Hypothesis ($\Hmuinfty$) holds,
the measure $\mu_\infty$ is {\em invariant} under the semigroup
$\PP$, that is, for all $f\in B_b(E)$ we have 
$$\int_E (P(t)f)(x)\,d\mu_\infty(x) = \int_E f(x)\,d\mu_\infty(x), \qquad t\ge
0.$$
By standard arguments, cf. \cite[Theorem XIII.1]{Yo}, 
it follows that $\PP$ extends to a 
$C_0$-contraction  semigroup, also denoted by $\PP$, on
$L^p(E,\mu_\infty)$ for all $p\in [1,\infty)$. 
The space $\F_0$, being norm-dense and $\PP$-invariant, is a core for the 
generator $(L,\D(L))$. 
We have the following integration by parts formula:

\begin{lemma}[$\Hmuinfty$]\label{part_int}  
For all $f,g\in \F_0$ we have
\begin{equation}\label{parts}
 \int_E f\, L g  + g \, L f \, d\mu_\infty  = -\int_E [D_H f, D_H
g]_{H} \,d\mu_\infty.
\end{equation}
\end{lemma}
\begin{proof}
Observe that $\F_0$ is closed under multiplication. Hence if $f,g\in \F_0$, 
then 
$fg\in \F_0\subseteq \D(L_0)\subseteq \D(L)$ 
and a simple caluculation based on \eqref{L0} 
gives
\begin{equation}
\label{eq:L(fg)}
\begin{aligned}
L(fg) =  L_0(fg) & =  f\, L_0 g  + g \, L_0 f  +  [D_{H} f, D_{H}
g]_{H}
 \\ & =  f\, L g  + g \, L f  +  [D_{H} f, D_{H} g]_{H}.
\end{aligned}
\end{equation}
Since $\mu_\infty$ is an invariant measure, we have 
\[\int_E P(t)(fg)\,d\mu_{\infty}=\int_E fg\,d\mu_{\infty}\]
from which it is immediate that 
\[\int_E L (fg)\,d\mu_{\infty}=0.\]
Therefore, for $f,g\in \F_0$ the desired result follows by integrating 
\eqref{eq:L(fg)} over $E$.
\end{proof}

\begin{remark} 
The identity \eqref{parts} extends to 
arbitrary elements $f,g\in \D(L)$ if $D_{H}$ is closable. 
Necessary and sufficient conditions for closability of $D_{H}$, as well as
simple examples where $D_{H}$ fails to be closable, were obtained 
in \cite{clos}. In Proposition \ref{closV} below we show that $D_H$
is closable if $\PP$ is analytic on $L^2(E,\mu_\infty)$. 
\end{remark}

On $L^2(E,\mu_\infty)$ we have the representation
$$P(t) = \Gamma (S_\infty\s(t)), \qquad t\ge 0,$$
where $\Gamma$ denotes the second quantization functor;
cf. \cite{CG1}, \cite{Ne}.
This result permits one to study the semigroup $\PP$ through the semigroup 
$\SS_\infty\s$. 
We give to simple illustrations.
The first is a characterization of selfadjointness.

\begin{theorem} [$\Hmuinfty$] The following assertions are
equivalent:

\begin{enumerate}
\item The semigroup $\PP$ is selfadjoint on $L^2(E,\mu_\infty)$;
\item The semigroup $\SS$ is $Q$-symmetric.
\end{enumerate}
\end{theorem}
\proof
We will show that $\SS$ is $Q$-symmetric if and only if
$\SS_\infty$ is selfadjoint. The proposition is then a consequence of
the identities
$P(t) = \Gamma(S_\infty\s)$ and $P\s(t)= \Gamma(S_\infty)$,
where $\Gamma$ denotes the second quantization functor.

If $\SS$ is $Q$-symmetric, then for all $t\ge 0$ and $x\s,y\s\in E\s$ we have
$$
\begin{aligned}
 {[S_\infty(t) \i\s x\s, \i\s y\s]_{\H}}
& = [ \i\s x\s, S_\infty\s(t)\i\s y\s]_{\H}
= [ \i\s x\s, \i\s S\s(t) y\s]_{\H}
\\ & = \lb \Q x\s, S\s(t) y\s\rb
= \lb \Q S\s(t)y\s, x\s\rb
\\ & = \lb S(t)\Q y\s, x\s\rb
 = [\i\s y\s, \i\s S\s(t)x\s]_{\H} 
\\ & = [\i\s y\s,  S_\infty\s(t)\i\s x\s]_{\H} 
= [S_\infty\s(t)\i\s x\s, \i\s y\s]_{\H}.
\end{aligned}
$$ 
It follows that $\SS_\infty$ is selfadjoint.
Conversely if $\SS_\infty$ is selfadjoint, then a similar argument
shows that $\SS$ is $Q$-symmetric.
\qed

The second illustration concerns the spectral gap of the generator of $\PP$.
If $(B,\D(B))$ is a negative operator in a Hilbert space $K$, i.e., if
$ [Bk,k]_K \le 0$ for all $k\in \D(B),$ 
then its spectrum $\sigma(B)$ is contained in the interval $(-\infty, 0]$.
We say that $B$ has a {\em spectral gap} if $0\in \sigma(B)$
and there exists $\omega>0$ 
such that $\sigma(B)\setminus\{0\}\subseteq (-\infty,-\omega]$.
The largest such $\omega>0$ is called the {\em spectral gap}
of $B$.

As an application of the results of Section \ref{sec:gap} we shall give a 
necessary and sufficient condition for 
the existence of a spectral gap for the generator $L$ of the Ornstein-Uhlenbeck
semigroup $\PP$ in $L^2(E,\mu_\infty)$.

Let ${\mathscr H}_1=L^2(E,\mu_\infty)\ominus\R$ denote the orthogonal complement 
in $L^2(E,\mu_\infty)$ of the constant functions. By second quantization
and the properties of the Wiener-It\^o decomposition, 
we obtain immediately that $\PP$ restricts to a $C_0$-semigroup 
of contractions on ${\mathscr H}_1 $ satisfying
$\|P(t)\|_{{\mathscr H}_1}=\|S_{\infty}(t)\|$ for all $t\ge 0$.
Let us denote the generator of $\PP_1$ by $L_1$. The following 
result may now be deduced 
from Lemma \ref{thm:gap}, Theorem \ref{thm:gap-2}, and the spectral
theory of $C_0-$semigroups. 

\begin{theorem}[$\Hmuinfty$] 
%
The following assertions are equivalent:
\begin{enumerate}
\item $L_1$ has a spectral gap; 
\item $\H\subseteq H$.
\end{enumerate}
If $H_{\infty}\subseteq H$, the spectral gap of $L_1$ equals
the exponential growth bound of the
semigroup generated by $A_\infty$.
\newline
If $\SS$ restricts to a $C_0$-semigroup on $H$, then $(1)$ and $(2)$ 
are equivalent to:
\begin{enumerate}
\item[(3)] $\SS_H$ is uniformly exponentially stable.  
\end{enumerate}
\end{theorem}

\section{Analyticity of the Ornstein-Uhlenbeck semigroup}
\label{sec:analytic}

In this section we investigate conditions under which the complexified semigroup
$\PP^{\C} = \{P^\C(t)\}_{t\ge 0} $
is analytic on $L_2^\C(E,\mu_\infty)$, the complexification of 
$L^2(E,\mu_\infty)$. Here $\PP$ denote the Ornstein-Uhlenbeck semigroup 
on $L^2(E,\mu_\infty)$ associated with $\SS$ and $H$; cf. Section
\ref{sec:preliminaries}.

Recall that a semigroup $\TT = \{T(t)\}_{t\ge 0}$ on a complex Banach
space is called an {\em analytic
contraction semigroup} if $\TT$ is analytic and 
$\n T(z)\n\le 1$ for all $z\in\C$ 
belonging to some sector containing the positive real axis.

Our first result generalizes to Banach spaces a result from \cite{go};
cf. also \cite[Theorem 3.6]{fu}.

\begin{theorem}[$\Hmuinfty$]\label{thm:analytic}
The following assertions are equivalent:
\begin{enumerate}
\item $\PP^\C$ extends to an analytic semigroup
on $L_2^\C(E,\mu_\infty)$;
\item $\PP^\C$ extends to an analytic contraction semigroup
on $L_2^\C(E,\mu_\infty)$;
\item $\SS_\infty^\C$ 
extends to an analytic contraction 
semigroup on 
$H_\infty^\C$;
\item There exists a constant $M \ge 0$ such that 
$$ \bigl | [A_\infty\s h_\infty , h_\infty']_{H_\infty} \bigr| 
\le M 
 \bigl | [A_\infty\s h_\infty, h_\infty]_{H_\infty} \bigr|^{\frac12} 
 \cdot \bigl | [A_\infty\s h_\infty', h_\infty']_{H_\infty}
 \bigr|^\frac{1}{2}
$$
for all $h_\infty, h_\infty'  \in \i\s \D(A\s)$.
\end{enumerate}
In this situation, $\PP^{\C}$ and $\SS_\infty^\C$ are contractive on the 
same sectors.
\end{theorem}

\proof
The proof is analogous to the corresponding result for Hilbert spaces $E$ given 
in \cite{go}, so we only sketch the main steps.

The equivalences
(1) $\,\Leftrightarrow\,$ (2) $\,\Leftrightarrow\,$ (3) as well as the 
final statement follow from the fact that $P^\C\!(t) =
\Gamma^\C\bigl((S_\infty^\C(t))\s\bigr)$, where $\Gamma^{\C}$ denotes the 
complex second quantization functor.

(3) $\, \Leftrightarrow\,$ (4): 
Since $ i_\infty\s(\D(A\s))$ 
is a core for $\D(A_\infty\s)$, the estimate in condition (4) 
holds for all $h_\infty,h_\infty'\in\D(A_\infty\s)$.  
 Hence by \cite[Proposition I.2.17]{MR}, 
condition (4) holds if and only if there exists $b>0$ such that
$$ [(A_\infty^\C)\s h^{\C}, h^{\C}]_{H_\infty^\C} 
\in \{z\in\C: \ |\Im z| \le b\, \Re z\}$$ 
for all
$h^{\C}\in\D((A_\infty^\C)\s).$
Since $1\in\varrho(A_\infty^\C)$ (recall that $\SS_\infty$, hence also $\SS_\infty^\C$, 
is a contraction semigroup),
by \cite[Theorem 1.5.9]{Go} this condition is in turn equivalent to
condition (3).
\qed

\begin{remark} Using the terminology of \cite{MR}, condition (4) says that $A_\infty$
satisfies a strong sector condition.
\end{remark}

We will develop Theorem \ref{thm:analytic} a little further.

\begin{theorem}[$\Hmuinfty$]\label{thm:holo-1} 
The following assertions are equivalent:

\begin{enumerate}
\item $\PP^\C $ extends 
to a analytic semigroup on $L_2^\C(E,\mu_\infty)$;
\item For all $x\s\in \D(A\s)$ we have  $AQ_\infty x\s\in H$, and  
there is a constant $C\ge 0$ such that 
$$\n AQ_\infty x\s\n_H \le C\n i\s x\s \n_H, \qquad x\s\in \D(A\s).$$ 
\item For all $x\s \in \D(A\s)$ we have 
$Q_\infty(A\s x\s )\in H$, and there exists a constant $C>0$ such that
$$\n Q_\infty A\s x\s\n_H \le C \n i\s x\s \n_H, \qquad x\s\in \D(A\s).$$
\end{enumerate}
\end{theorem}
\proof
We start by noting that for all $x\s\in \D(A\s)$ and $y\s\in \D(A\s)$
we have $i_\infty\s x\s\in\D(A_\infty\s)$, $i_\infty\s y\s\in\D(A_\infty\s)$,
$$\lb A Q_\infty x\s, y\s\rb = \lb \Q A\s y\s,x\s\rb = [i_\infty\s A\s y\s,
i_\infty\s x\s]_{\H} = [A_\infty\s i_\infty\s y\s,
i_\infty\s x\s]_{\H},$$
and
$$\lb Qx\s,x\s\rb = -2 \lb \Q A\s x\s,x\s\rb
=-2[A_\infty\s i_\infty\s x\s, i_\infty\s x\s]_{H_\infty} .$$
Hence by Theorem \ref{thm:analytic}, $\PP^\C$ is analytic if and only if
there exists a constant $M\ge 0$ such that
\begin{equation}
\label{reformulation}
|\lb A Q_\infty x\s, y\s\rb |\le M \lb  Qx\s, x\s \rb^{\frac12} \lb  Qy\s, y\s\rb^{\frac12}, \qquad 
x\s, y\s\in \D(A\s).
\end{equation}

(1)$\,\Rightarrow\,$(2):
By Theorem \ref{thm:analytic}, for all $x\s, y\s\in \D(A\s)$ we have
\begin{equation}\label{eq:kernels}
|\lb AQ_\infty x\s, y\s \rb| \le  M\n i\s x\s\n_H\, \n i\s y\s\n_H.
\end{equation}
It follows that the map
$i\s y\s \mapsto \lb  AQ_\infty x\s, y\s\rb$
is well defined and can be extended to a bounded linear form on $H$ of norm $\le \n  i\s x\s\n_H$. 
Therefore by the Riesz representation theorem we can identify
$AQ_\infty x\s$ with an element of $H$ of norm $\le \n i\s x\s\n_H$.
This gives (2), with $C=M$.

(2)$\,\Rightarrow\,$(1):
For all $x\s, y\s\in \D(A\s)$ we have
$$
\begin{aligned}  
 |\lb AQ_\infty x\s, y\s\rb| & = |[AQ_\infty x\s, i\s y\s]_H| 
 \\ & \le C \n i\s x\s\n_H \n i\s y\s\n_H 
 = C\lb Qx\s,x\s\rb^{\frac12}\lb Qy\s,y\s\rb^{\frac12}
\end{aligned}
$$
and  $\PP^\C $ has an analytic extension.

(1)$\,\Rightarrow\,$(3): By the proof of (1)$\,\Rightarrow\,$(2),
for all $x\s\in \D(A\s)$ we have  
$Q_\infty A\s x\s = - Qx\s - AQ_\infty\s x\s\in H$,  
and there is a  constant $c$ such that
$$ \n Q_\infty A\s x\s\n_H  \le  \n Qx\s \n_H + \n A\Q x\s \n_H \le (1+c)\n
i\s x\s\n_H$$
for all $x\s\in \D(A\s)$.

(3)$\,\Rightarrow\,$(1): 
For all $x\s, y\s\in \D(A\s)$ we have
$$
\begin{aligned}  
|\lb A \Q x\s, y\s\rb |=| \lb \Q A\s y\s, x\s\rb | 
& = |[\Q A\s y\s, i\s x\s]_H| 
\\ &  \le C \n i\s x\s\n_H \n i\s y\s\n_H
= C\lb Qx\s,x\s\rb^{\frac12}\lb Qy\s,y\s\rb^{\frac12}.
\end{aligned}
$$
\qed

This leads to the following concise {\em necessary} condition for
analyticity of the Ornstein-Uhlenbeck semigroup:

\begin{corollary}[$\Hmuinfty$]\label{necessary}
If $\PP^\C $ extends to an analytic semigroup on $L_2^\C(E,\mu_\infty)$, 
then  for all $x\s\in \D(A\s)$ with $Qx\s = 0$ we have $QA\s x\s = 0$.
\end{corollary}
\proof
If $Qx\s = 0$, then $i\s x\s = 0$ and hence $\Q A\s x\s = 0$. 
This we combine with the simple observation that $\ker \Q \subseteq \ker
Q$, cf. \cite[Lemma 5.2]{clos}.
\qed

\begin{example}
Let $E=\R^2$ and let $Q$ and $\SS$ be as in Example \ref{ex:Hinfty}.
Since Hypothesis $(\HQinfty$) holds and $E$ is finite-dimensional,
Hypothesis $(\Hmuinfty$) trivially holds as well. 
Denote the centred Gaussian measure associated with $Q_\infty$ by
$\mu_\infty$.
By direct computations, Fuhrman \cite{fu} showed that the 
associated Ornstein-Uhlenbeck semigroup $\PP^\C$ fails to be analytic
on $L_2^\C(E,\mu_\infty)$. We derive this from Corollary \ref{necessary},
by noting that
$$ Q \left(\begin{array}{c}1 \\ 0 \end{array}\right)
   = \left(\begin{array}{c}0 \\ 0\end{array}\right),
   \qquad 
   QA\s \left(\begin{array}{c}1 \\ 0 \end{array}\right) 
   = Q\left(\begin{array}{r}-1 \\ 1\end{array}\right) 
   =  \left(\begin{array}{c}0 \\     1  \end{array}\right).
$$
Notice that $\SS_\infty$ is both contractive 
and analytic (its generator being bounded).
Hence the 
same is true for its complexification $\SS_\infty^\C$. This does not
contradict Theorem \ref{thm:analytic}; the point is that
$\SS_\infty^\C$ fails to be an analytic contraction semigroup in the sense 
of the definition given at the beginning of this section.
This can be verified explicitly by extending the computation in Example
\ref{ex:Hinfty} to complex time. By doing so we obtain
$$
 \n S_\infty^\C(z)\n = e^{-\Re z } \bigl(|z| + \sqrt{|z|^2+1}\bigr).
$$
Let $z=re^{i\theta}$ for a certain $\theta\in
\left(-\frac{\pi}{2},\frac{\pi}{2}\right)\setminus\{0\}$. 
Then 
$$ \n S_\infty^\C(z)\n = e^{-r\cos\theta}\left(r+\sqrt{r^2+1}\right).$$
We claim that for any $\theta\in\left(0,\frac{\pi}{2}\right)$ 
we have 
$$e^{-r\cos\theta}\left(r+\sqrt{r^2+1}\right)>1\eqno (4.1)$$
for all sufficiently small $r>0$.
Indeed, (4.1) holds if and only if 
$$f(r)=r+\sqrt{r^2+1}>e^{r\cos\theta}=: g_\theta(r)$$
for some $r>0$. This is clearly true for  small  $r>0$ 
because $f(0)=g_\theta(0)=1$ 
and $f'(0)=1>g_\theta'(0)=\cos\theta$. 
\end{example}

In the next corollary, which is a minor extension of
 \cite[Corollary 2.5]{go}, we specialise Theorem \ref{thm:holo-1} 
to Hilbert spaces $E$. We identify $E$ and its dual in the usual way.

\begin{corollary}[$\Hmuinfty$]
Suppose $E$ is a Hilbert space and $Q\in {\mathscr L}(E)$ has a bounded inverse.
Then the following assertions are equivalent:
\begin{enumerate}
\item 
The semigroup $\PP^\C $ extends to a analytic semigroup on 
$L_2^\C(E,\mu_\infty)$;
\item The operator $AQ_\infty$ extends to a bounded operator on $E$;
\item The operator $Q_\infty A\s$ extends to a bounded operator on $E$.
\end{enumerate}
\end{corollary}

The final result of this section is closely related to \cite[Proposition 3.3]{MR}.

\begin{proposition}[$\Hmuinfty$]\label{closV}
If the transition semigroup $\PP^\C$ is analytic, then $D_H$ is closable.
\end{proposition}
\begin{proof}
We introduce a
densely defined operator $(V,\D(V))$ from $\H$ to $H$,
$$
\begin{aligned}
\D(V)       & := \{\i\s x\s: \ x\s\in E\s\}, \\
  V(\i\s x\s) & := i\s x\s, \qquad x\s\in E\s.
\end{aligned}
$$
It was shown in \cite{clos} $D_H$ is closable in
$L^2(E,\mu_\infty)$ if and only $V$ is closable.

We will show that $V$ is closable if $\PP^\C$ is analytic. The 
proof uses the trick from \cite[Theorem 2.15]{MR}.
Let $\i\s x_n\s\to 0$ in $\H$ and $V (\i\s x_n\s) = i\s x_n\s \to g$ 
in $H$; we have to prove
that $g=0$. Fix $\e>0$ arbitrary and 
choose an index $N$ large enough such that
$\n V(\i\s x_n\s - \i\s x_N\s)\n_H\le \e$ and
$\n \i\s x_n\s\n\le \e$ for all $n\ge N$.
Then for all $n\ge N$ we have
$$
\begin{aligned}
 \tfrac12\n V(\i\s x_n\s) \n_H^2 & =\tfrac12 \lb Qx_n\s, x_n\s\rb
 = |\lb \Q A\s x_n\s, x_n\s\rb|
\\ & \le   |\lb \Q A\s (x_n\s-x_N\s), (x_n\s-x_N\s)\rb|
\\ &\qquad    + |\lb \Q A\s (x_n\s-x_N\s), x_N\s\rb|
    + |\lb \Q A\s x_N\s, x_n\s\rb|
\\ &= \tfrac12 \n V(\i\s x_n\s-\i\s x_N\s)\n^2
\\ &\qquad + |[\Q A\s (x_n\s-x_N\s), i\s x_N\s]_{H}|
  + |[ \i\s A\s x_N\s, \i\s x_n\s]_{\H}|
\\ & \le  \tfrac12 \n V(\i\s x_n\s-\i\s x_N\s)\n^2 
\\ &\qquad   + C\n i\s (x_n\s-x_N\s)\n_{H}  \n i\s x_N\s\n_{H}
  + \n \i\s A\s x_N\s\n_{\H} \n \i\s x_n\s\n_{\H}
\\ & =  \tfrac12 \n V(\i\s x_n\s-\i\s x_N\s)\n^2 
\\ &\qquad   + C\n V(\i\s x_n\s-\i\s x_N\s)\n_{H}  \n i\s x_N\s\n_{H}
  + \n \i\s A\s x_N\s\n_{\H} \n \i\s x_n\s\n_{\H}
\\ & \le \tfrac12\e^2 + C\e M + \n \i\s A\s x_N\s\n_{\H} \n \i\s
x_n\s\n_{\H},
\end{aligned}
$$
where $C$ is the constant from Theorem \ref{thm:holo-1}(3) and 
$M:= \sup_n \n i\s x_n\s\n_H$ is finite since $\limn i\s x_n\s =
g$. 
Upon letting
$n\to\infty$, it follows that
$$ \limsup_{n\to\infty} \ \tfrac12\n V(\i\s x_n\s) \n_H^2 
\le \tfrac12\e^2 + CM\e.
$$
Since $\e>0$ was arbitrary, we conclude that
$ g = \limn V(\i\s x_n\s) =0.$ 
\end{proof}

\section{Analyticity and invariance of $H$}
\label{sec:analytic_H1}

It turns
out that there is a close relationship between analyticity of the 
Ornstein-Uhlenbeck semigroup and invariance of $H$. This will be the 
topic of the present section.

We start with a necessary condition for analyticity:

\begin{theorem}[$\Hmuinfty$]\label{thm:analyticH_2}
If $\H\subseteq H$ and $\PP^\C$ is analytic, then $\SS^\C$ restricts to a 
bounded analytic $C_0$-semigroup on $H^\C$.
\end{theorem}
\begin{proof}\footnote{ \ The proof in the published version of this paper contains
a mistake which has been corrected here.} 
By Proposition \ref{Riesz} 
it suffices to check that there exists a constant $M$ such that 
for all $t\ge 0$ and $x\s\in E\s$ we have
$$ \n i\s S\s(t) x\s\n_H \le M \n i\s x\s\n_H.$$
Since $A_\infty\s$ generates an analytic $C_0$-contraction semigroup on
$H_\infty$, the form
$$ \mathscr{E}(g,h) := -[A_\infty\s g,h]_{H_\infty}, \quad
g,h\in\D(A_\infty\s),$$ is sectorial. 

By Proposition \ref{closV} the operator $V$ 
is closable as a densely defined operator from $H_\infty$ to $H$.
Let $\overline{V}$ be its closure. 
Its domain, $\D(\overline{V})$,  
 is a Banach space with respect to the graph norm
$ \n h\n_{\D(\overline{V})}^2 := \n h\n_{H_\infty}^2 + \n
\overline{V}h\n_{H}^2$.
Taking $h = i_\infty\s x\s$ and using \eqref{eq:incl} 
we have two-sided estimate
$$
\begin{aligned} 
\n \overline{V}i_\infty\s x\s\n_{H}^2
  \le  \n i_\infty\s x\s \n_{\D(\overline{V})}^2 
 & = \n i_\infty\s x\s \n_{H_\infty}^2 + \n \overline{V}i_\infty\s x\s\n_{H}^2
\\ & \le K^2 \n i\s x\s \n_{H}^2 + \n \overline{V}i_\infty\s
x\s\n_{H}^2 
= (K^2+1) \n \overline{V}i_\infty\s x\s \n_{H}^2,
\end{aligned}
$$
where $K$ is the norm of the embedding $H_\infty\embed H$.
This shows that 
$$ \nn h \nn_{\D(\overline{V})} := \n \overline{V}h\n_{H}, \quad
h\in\D(\overline{V}), $$
defines an equivalent norm on  $\D(\overline{V})$.

We claim that $\D(\overline{V})$ can be identified with the form domain of 
${\mathscr E}$. By Theorem \ref{thm:holo-1} the mapping
$B: i\s x\s \mapsto Q_\infty A\s x\s$, defined on the dense subspace $i\s
\D(A\s)$ of $H$, takes values in $H$ and extends to a bounded operator 
$B$ on $H$. Moreover, for all 
$x\s, y\s\in\D(A\s)$ we have $i_\infty\s x\s \in \D(A_\infty\s)$ and 
$$\mathscr{E}(i_\infty\s x\s,i_\infty\s y\s)   = 
-[A_\infty\s i_\infty\s x\s, i_\infty\s y\s]_{H_\infty} = -[B{V}i_\infty\s x\s, 
{V}i_\infty\s y\s]_H.$$
Therefore,
for all $h\in \D(\overline{V})$, 
$$\mathscr{E}(g,h)   = -[B\overline{V}g, \overline{V}h]_H, \quad g,h\in \D(\overline{V}).$$
This proves the claim.

It follows from the general theory of sectorial
operators that $\D(\overline{V})$ is invariant under $S_\infty\s$ and that the
restriction of $S_\infty\s$ to $\D(\overline{V})$ is a bounded 
analytic $C_0$-semigroup. 
Therefore, for some constant $m$ and  all $t\ge 0$ and $x\s\in E\s$,
$$
\begin{aligned}
\n i\s S\s(t) x\s\n_H & =
\n \overline{V}i_\infty\s S\s(t)x\s\n_{H} 
 = \n \overline{V}S_\infty\s (t)i_\infty\s x\s\n_{H} 
\\ & = \nn S_\infty\s (t) i_\infty\s x\s\nn_{\D(\overline{V})}
 \le  c\nn i_\infty\s x\s\nn_{\D(\overline{V})}
 = c\n \overline{V}i_\infty\s x\s\n_{H}= c\n i\s x\s\n_H.
\end{aligned}
$$
This proves that $H$ is $S$-invariant and that the restricted semigroup $S_H$ is 
bounded. By Proposition \ref{prop:str-cont}, $S_H$  is strongly continuous.
It remains to prove that $S_H$ is bounded analytic on
$H$.
Since $S_\infty\s$ restricts to a bounded analytic $C_0$-semigroup on
$ \D(\overline{V})$ there is a constant $C$ such that for all $t>0$ and
$h\in\D(\overline{V})$,
$$ \n  \overline{V} A_\infty\s S_\infty\s(t)h\n_H \le \frac{C}{t}\n  \overline{V}h\n_H.$$
As above, taking $h = i\s x\s$ this implies 
$$ \n i\s A\s S\s(t) x\s\n_H \le  \frac{C}{t} \n i\s x\s\n_H.$$
Therefore, $\n A_H\s S_H\s(t)\n_H \le  \frac{C}{t}$, which implies the result.
\end{proof}

We proceed with a partial converse.
\begin{theorem}[$\Hmuinfty$]\label{thm:analyticH_1}\footnote{ \ In the published version of the paper, 
the word
`analytic' was missing in the first line of the statement of the result.}
Suppose that $\SS$ restricts to an analytic $C_0$-semi\-group $\SS_H$ on $H$
which is contractive in some 
equivalent Hilbertian norm on $H^{\C}$.
Then $\PP^\C$ is analytic.
\end{theorem}

\begin{proof}
We will show that 
$\SS_\infty^\C$ is an analytic contraction semigroup on $H_\infty^\C$.
Once this is proved, the theorem follows by an appeal to Theorem
\ref{thm:analytic}.

Identifying $H$ and its dual in the usual way, we 
define $R_\infty\in \L(H)$ by
$$ R_\infty h := \int^\infty_0 S_H(t)S_H\s h\, dt\qquad (h\in H).$$
Let $i:H\embed E$ denote the embedding; then we have $\Q = i\circ R_\infty \circ
i\s$. By an observation in Section \ref{sec:preliminaries}, 
the RKHS's $\H = H_{\Q}$ and
$H_{R_\infty}$ are canonically isometrically isomorphic as Hilbert spaces
and identical as subsets of $E$. By complexifying, the same is true 
for their complexifications $H_\infty^{\C}$ and $H_{R_\infty}^{\C}$.
It follows that, in order to prove that $\SS_\infty^\C$
extends to an analytic contraction semigroup on $H_\infty^{\C}$, 
{\em we may assume without loss of generality that $E = H$ and $Q = I$}.
We also note that $\H = H_{R_\infty}\subseteq H$.

Let $\nn \cdot\nn$ be an equivalent Hilbertian norm on $H^{\C}$
such that 
$ \nn S_H^{\C}(z)\nn \le 1$
for all $z$ in some sector containing the positive real axis, and let
$[\![\cdot,\cdot]\!]$ be the corresponding inner product.
For all $x,y\in H$ we have
$$ \nn x+iy\nn^2 = [\![x+iy, x+iy]\!] = [\![x, x]\!] + [\![y,y]\!] = 
\nn x\nn^2 + \nn y\nn^2.
$$
Hence $(H^{\C}, \nn\cdot\nn)$ is the complexification of its real part, 
and we may apply the observation from Section \ref{sec:preliminaries} 
once more, this time to the isomorphism $j: (H, \n\cdot\n) \simeq (H, \nn\cdot\nn)$.
It follows that the RKHS's associated with $R_\infty$ and $j\circ R_\infty\circ
j\s$ are canonically isometrically isomorphic, and identical as subsets 
of $H$, and again the same is true for their complexifications. 
Thus, in order to prove that $\SS_\infty^\C$
extends to an analytic contraction semigroup on $H_\infty^{\C}$, 
{\em it even suffices to prove this for the case where $\SS_H^{\C}$ extends
to an analytic {\em contraction} semigroup on $H^\C$}. 

It is well known that
$$ H_\infty = \left\{\int^\infty_0 S_H(t)f(t)\,dt: \ f\in L^2(\R_{\!+};H)  \right\} $$
with norm given by
$$ 
 \n h_\infty\n_{\H} \ = \ \inf\left\{\n f\n_{L^2(\R_{\!+};H)}: 
 \  h_\infty = \int^\infty_0 S_H(t)f(t)\,dt \right\},
$$
cf. \cite[Appendix B]{DZ}.
Upon complexifying we see that
$$
 H_\infty^{\C} = \left\{\int^\infty_0 S_H^{\C}(t)f(t)\,dt: \ f\in
 L^2(\R_{\!+};H^{\C})  \right\} 
$$
with norm given by
\begin{equation}\label{normC}
 \n h_\infty\n_{H_\infty^{\C}} \ = \ \inf\left\{\n f\n_{L^2(\R_{\!+};H^{\C})}: 
 \  h_\infty = \int^\infty_0 S_H^{\C}(t)f(t)\,dt \right\}.
\end{equation}
Indeed, the representation of $H_\infty^{\C}$ follows immediately by
considering real and imaginary parts of elements in $H_\infty^\C$ separately, 
and the expression \eqref{normC} for the
complexified norm is proved as follows. Denote the infimum on the right hand
side of \eqref{normC} by $I^{\C}$.
Fix $h_\infty \in H_\infty^\C$ and write $h_\infty = a_\infty + i b_\infty$ with
$a_\infty, b_\infty\in \H$. Fix $\e>0$ arbitrary 
and choose $f,g\in L^2(\R_{\!+},H)$ representing $a_\infty, b_\infty$ such that
$$
   \n f\n_{L^2(\R_{\!+},H)}^2\ 
  \le \ \n a_\infty\n_{\H}^2 + \e, \qquad
  \n g\n_{L^2(\R_{\!+},H)}^2\ 
  \le \ \n
 b_\infty\n_{\H}^2 + \e.
$$
Then,
$$
\begin{aligned}
 \n f+ig \n_{L^2(\R_{\!+},H^{\C})}^2
  & = \n f\n_{L^2(\R_{\!+},H)}^2 + \n g\n_{L^2(\R_{\!+},H)}^2 
 \\ & \le (\n a_\infty\n_{\H}^2+\e) + (\n b_\infty\n_{\H}^2+ \e) 
 = \n h_\infty\n_{H_\infty^\C}^2+2\e,
\end{aligned}
$$
which shows that
$
I^{\C}\le  \n h_\infty\n_{H_\infty^\C}.
$
On the other hand, if $f,g\in L^2(\R_{\!+},H)$ are arbitrary functions
representing $a_\infty, b_\infty$, then
$$ \n h_\infty\n_{H_\infty^{\C}}^2
 =   \n a_\infty\n_{\H}^2 + \n b_\infty\n_{\H}^2 
 \le \n f\n_{L^2(\R_{\!+},H)}^2 + \n g\n_{L^2(\R_{\!+},H)}^2
 = \n f+ig \n_{L^2(\R_{\!+},H^{\C})}^2,
$$ 
which gives the converse inequality
$\n h_\infty\n_{H_\infty^{\C}} \le I^{\C}.$ 
This proves \eqref{normC}.

Now it is easy to finish the proof. Given $h_\infty\in H_\infty^{\C}$,
choose  an arbitrary $f\in L^2(\R_{\!+};H^\C)$ representing $h_\infty$:
$$ h_\infty = \int^\infty_0 S_H^\C(t)f(t)\,dt.$$
Then, for any $z$ in the sector where $\SS_H^{\C}$ is contractive, we have
$$ S_H^\C(z)h_\infty = S_H^\C(z)\int^\infty_0 S_H^\C(t)f(t)\,dt
= \int^\infty_0 S_H^\C(t)[ S_H^{\C}(z)f(t)]\,dt.$$
It follows that $S_H^\C(z)h_\infty\in H_\infty^{\C}$, with norm
$$ \n S_H^\C(z)h_\infty\n_{H_\infty^{\C}} \le \n
S_H^{\C}(z)f(\cdot)\n_{L^2(\R_{\!+};H^\C)} \le \n f(\cdot)\n_{L^2(\R_{\!+};H^\C)}.
$$
Taking the infimum over all representing functions $f$, we obtain
$$  \n S_H^\C(z)h_\infty\n_{H_\infty^{\C}} \le\n h_\infty\n_{H_\infty^{\C}} .$$
It follows that the operators $S_H^\C(z)$ 
restrict to a contractions on $H_\infty^\C$.
The restriction of $\SS_H^\C$ to $H_\infty^\C$
agrees with $\SS_\infty^\C$ for real time,
and it is routine to check that 
it is strongly continuous and analytic.
\end{proof}

Notice that there is only a small gap between Theorems 
\ref{thm:analyticH_2} and \ref{thm:analyticH_1}. The
assumption $\H\subseteq H$ in Theorem 
\ref{thm:analyticH_2} implies that 
$\SS_H$ is uniformly exponentially stable, and conversely
the assumption in Theorem 
\ref{thm:analyticH_1} that
$\SS_H$ is uniformly exponentially stable implies that $\H\subseteq H$.

The assumptions of Theorem \ref{thm:analyticH_1} are fulfilled when
$E$ is a Hilbert space, $H=E$, and $\SS$
satisfies an estimate of the type
$\n S(t)\n\le e^{-\omega t}$ for some $\omega>0$ and all $t\ge 0$. 
In this special setting, 
the theorem is due to Da Prato \cite{DP}, who proved it 
by using interpolation theory
and maximal regularity.

\begin{remark}
The existence of an equivalent Hilbertian norm on $H^\C$
in Theorem \ref{thm:analyticH_1}
is equivalent to the existence of an isomorphism $T:H^\C\to H^\C$
such that 
\begin{equation}\label{isom}
\n T^{-1}S_H^\C(z)T\n\le 1
\end{equation}
for al $z$ in some sector containing the positive real line.
The question when such an isomorphism exists is related to a 
famous question posed by Halmos in \cite{halmos}. For bounded 
analytic semigroups, this question was answered recently by Le Merdy 
\cite{lemerdy}. To quote his answer let us recall first that if 
${\rm ker}\,A=\{0\}$ and $A$ generates a bounded analytic semigroup 
on $H^\C$, then for any $s\in{\mathbb R}$ one can 
define a closed operator $(-A)^{is}$. We say that $A$ has {\em 
bounded imaginary powers} (briefly, $A\in {\rm BIP}$) 
if $(-A)^{is}$ is bounded for all $s\in\R$ and the function
$s\mapsto \n (-A)^{is}\n$ is locally bounded on $\R$.

It is known that $A\in{\rm BIP}$ in the following important cases:
\begin{itemize}
\item[(a)] $A$ is $m$-dissipative on $H$; 
\item[(b)] $A$ is normal and sectorial on $H$;
\item[(c)] $A$ generates a bounded $C_0$-group on $H$.
\end{itemize}
By an example of Baillon and Cl\'ement \cite{BC}, there exist analytic semigroups
on Hilbert spaces which are uniformly bounded on a sector, but whose
generator does not belong to BIP.  

Le Merdy  \cite{lemerdy} proved that for a bounded analytic semigroup 
whose generator $A$ satisfies ${\rm ker}\,A\,=\{0\}$,  \eqref{isom} 
holds if and only if $A\in{\rm BIP}$. 

In the situation of Theorem \ref{thm:analyticH_1}, $\SS_H$ is uniformly 
exponentially stable and therefore $0\not\in\sigma(A_H)$.
Hence the condition  ${\rm ker}\,A_H=\{0\}$
is trivially fulfilled.
\end{remark}

In the example below we consider a stochastic linear heat equation with
correlated cylindrical noise in $L^p(\O)$ with $p\in [2,\infty)$.
Similar equations were considered in \cite{DaPr},
where it is a starting point for the analysis of nonlinear stochastic
differential equations with dissipative drifts.

\begin{example}
Let $\O$ be a bounded open domain in $\R\!{}^d$ with $C^2$-boundary,
let $2\le p<\infty$,
and let $A$ be the $L^p(\O)$-realization 
of a uniformly elliptic differential operator of the form 
$$
A_0 \, = \, \sum_{i,j=1}^d a_{ij}\,\partial_{ij} +
\sum_{i=1}^d b_i \,\partial_i
$$
with domain $\D(A)=W^{2,p}(\O)\cap W_0^{1,p}(\O)$. 
We assume 
that the coefficients $a_{ij} = a_{ji}$ belong to 
$C^{\theta}(\overline{\O})$ for a certain $\theta\in (0,1)$
and that the functions $b_i$ are bounded and 
measurable on $\O$. Under 
these assumptions it is known that $A$ generates a uniformly exponentially
stable and  analytic $C_0$-semigroup 
${\bf S}$ in $L^p(\O)$; cf. \cite{agmon}. 

In $E = L^p(\O)$ we consider a stochastic evolution equation
$$ dX(t) = A X(t) \,dt + dW_H(t). $$
Here $H$ is a separable Hilbert space which is continuously 
embedded into $E$ and $\{W_H(t)\}_{t\ge 0}$ is a (possibly cylindrical) 
Wiener process with Cameron-Martin space $H$.

We will consider two cases. To simplify notations, we will not 
distinguish between real spaces and their complexifications.

\medskip
%
%
\medskip
(a)  \ We take $p=2$, $E=L^2(\O)$, and $H= H_0^\beta(\O)$ with 
$\beta\ge 0$ and $\beta > \frac{d}{4} - \frac12$.
By a result 
in \cite{hiroshima} $A\in\mbox{\rm BIP}$. In case $\beta=0$ we can 
apply Le Merdy's result 
to find an equivalent Hilbertian norm in which $\SS_H = \SS$ is an analytic
contraction semigroup.
In case $\beta>0$, the fact that $A \in\mbox{\rm BIP}$ implies
that $H$ equals the interpolation space
$D_{A}\bigl(\frac {\beta}2,2\bigr)$ up to an equivalent norm.
Then by interpolation theory, 
${\bf S}_H$ is an analytic $C_0$-semigroup on $H$ which is 
contractive with respect to the
$D_{A}\bigl(\frac {\beta}2,2\bigr)$ norm. 
From \cite{agmon} we have 
\[\int_0^{\infty}\left\|S(t)\circ i_\beta\right\|_{\L_2(H_0^\beta(\O),L^2(\O))}^2\, dt
= \int_0^{\infty}\left\|A^{-\frac{\beta}{2}} S(t)\right\|_{\L_2(L^2(\O))}^2\, dt
 < \infty ,\]
where $\n \cdot\n_{\L_2}$ 
denotes the Hilbert-Schmidt norm
and $i_\beta: H_0^\beta(\O)\embed L^2(\O)$ is the inclusion mapping.
This implies that Hypothesis ($\Hmuinfty$) is satisfied.

By Theorem \ref{thm:analyticH_1}, the
associated Ornstein-Uhlenbeck  
semigroup ${\bf P}$ is analytic in $L^2(E,\mu_\infty)$. 

\medskip
(b) \ Let $p\in (2,\infty)$, $E = L^p(\O)$, and 
$H = H_0^{\alpha}(\O)$ with 
$\alpha > \frac12 - \frac1p$ (in dimension $d=1$)
or 
$\alpha> d(\frac34 - \frac1p) - \frac12$ (in dimensions $d\ge 2$). 
In both cases we may choose $\beta\ge 0$ with 
$\beta> \frac{d}{4}-\frac12$  and $\gamma> d(\frac12-\frac1p)$ 
such that
$\alpha> \beta+\gamma$.
By the Sobolev embedding theorem
we have a continuous inclusions $H\embed H_0^\gamma(\O)\embed E$. 
Then $H$ equals the interpolation space
$D_{A}\bigl(\frac {\alpha}2,2\bigr)$ up to an equivalent norm.  
Again  $H$ is invariant for ${\bf S}$, and
${\bf S}_H$ is an analytic $C_0$-semigroup on $H$ which is 
contractive with respect to the
$D_{A}\bigl(\frac {\alpha}2,2\bigr)$ norm. 

We will show next that Hypothesis ($\Hmu$) is satisfied in $E$.
The argument will be somewhat informal but can easily be rewritten 
in a rigorous way.
First recall that the realization of $A$ in $L^2(\O)$ belongs to BIP, 
from which it follows that $H = \D((-A)^{\frac{\a}{2}})$ with equivalent norms.
Hence, $(-A)^{\frac12(\a-\beta)}$ is an isomorphism from $H$ onto $H_0^\beta(\O)$,
and
$$W_{A}(t) := (-A)^{\frac12(\a-\beta)} W_H(t)$$
defines a cylindrical Wiener process whose Cameron-Martin
space equals $H_0^{\beta}(\O)$ up to an equivalent norm.
Then by case (b), the $L^2(\O)$-valued process  
$$ Y(t) := \int^t_0 S(t-s)\,dW_{A}(s), \qquad t\ge 0$$
solves the equation
$$ dY(t) = A Y(t)\,dt + dW_{A}(t) $$
with initial condition $Y(0)=0$ in $L^2(\O)$.
Then the process $\{X(t)\}_{t\ge 0}$ defined by 
$$X(t) := (-A)^{-\frac12(\a-\beta)}Y(t)$$ takes values in
$H$, hence in $E$, and solves the original
equation in $E$,
$$ dX(t) = A X(t)\,dt + dW_{H}(t) $$
with initial condition $X(0)=0$. It follows from Proposition \ref{prop:sACP} 
that Hypothesis ($\Hmu$) is satisfied in $E$. This proves the claim.
By \cite{NW},  the uniform exponential stability of $\SS$ in $E$ 
now implies that also ($\Hmuinfty$) is satisfied in $E$.

In conclusion, Theorem \ref{thm:analyticH_1} applies and we find that 
the Ornstein-Uhlenbeck semigroup 
${\bf P}$ is analytic in $L^2(E,\mu_\infty)$. 
\end{example}

\setcounter{footnote}{4}
{}\footnotetext{The references have been updated.}


\begin{thebibliography}{99}

\smallskip\bibitem{agmon} {\sc S.A. Agmon}, ``Lectures on Elliptic Boundary Value Problems'',  Van Nostrand, 1965.

\smallskip\bibitem{AN} {\sc W. Arendt and N. Nikolski}, 
{\em Vector-valued holomorphic functions revisited}, 
Math. Z. {\bf 234} (2000), 777--805.
 
\smallskip\bibitem{BC} {\sc J.-B. Baillon and Ph. Cl\'ement},
{\em Examples of unbounded imaginary powers of operators},
J. Funct. Anal. {\bf 100} (1991),  419--434.

\smallskip\bibitem{BRS} {\sc V. I. Bogachev, M. Rockner and B. Schmuland},
{\em Generalized Mehler semigroups and applications}, Probab. Theory Related Fields {\bf 105} (1996), 193--225

\smallskip\bibitem{Brz} {\sc Z. Brze\'zniak},
{\em Stochastic partial differential equations in M-type $2$ Banach spaces}, Potential Anal. {\bf 4} (1995),
1--45.

\smallskip\bibitem{BN} {\sc Z. Brze\'zniak and J.M.A.M. van Neerven},
{\em Stochastic convolution in separable Banach spaces
and the stochastic linear Cauchy problem}, 
Studia Math. {\bf 143} (2000), 43--74. 

\smallskip\bibitem{CG} {\sc S. Cerrai and F. Gozzi}, 
{\em Strong solutions of Cauchy problems associated to weakly
continuous semigroups}, 
Diff. Integral Eq. {\bf 8} (1995), 465--486. 

\smallskip\bibitem{CG1} {\sc A. Chojnowska-Michalik and B. Goldys}, 
{\em Nonsymmetric Ornstein-Uhlenbeck semigroup as second quantized operator}, 
J. Math. Kyoto Univ. {\bf 36} (1996), 481--498.

\smallskip\bibitem{CG2} {\sc A. Chojnowska-Michalik and B. Goldys},
{\em Nonsymmetric Ornstein-Uhlenbeck generators}, 
in: ``Infinite 
Dimensional Stochastic Analysis'' (Amsterdam, 1999) 
Verh. Afd. Natuurkd., 1st Series, Vol. 52, 
Royal. Netherl. Acad. Arts Sci., Amsterdam, 2000, pp. 99--116.

\smallskip\bibitem{CG3} {\sc A. Chojnowska-Michalik and B. Goldys}, 
{\em Symmetric Ornstein-Uhlenbeck semigroups and their generators},
Probab. Theory Relat. Fields. {\bf 124} (2002), 459--486.

\smallskip\bibitem{CW} {\sc K.L. Chung and R.J. Williams}, ``Introduction to Stochastic Integration'', 2nd ed., Birkh\"auser Verlag, 1990. 

\smallskip\bibitem{DP} {\sc G. Da Prato}, 
{\em Null controllability and strong Feller property 
of Markov transition semigroups}, 
Nonlinear Anal. {\bf 25} (1995), 941-949.

\smallskip\bibitem{DP2} {\sc G. Da Prato},  
{\em Bounded perturbations of Ornstein-Uhlenbeck semigroups},  
in: ``Evolution Equations, Semigroups and Functional Analysis'' (Milano, 2000),  97--114,
Progr. Nonlinear Differential Equations Appl., Vol. 50,
Birkh\"auser, Basel, 2002. 

\smallskip\bibitem{DaPr} {\sc G. Da Prato},  
``Stochastic Evolution Equations by Semigroup Methods'', 
Centre de Recerca Matem\`atica, Quaderns n\`um. 11, 1998.
 
\smallskip\bibitem{DZ} {\sc G. Da Prato and J. Zabczyk}, 
``Stochastic Equations in Infinite Dimensions'', 
Encyclopedia of Mathematics and its Applications, 
Cambridge University Press, Cambridge, 1992.

\smallskip\bibitem{DU} {\sc J. Diestel and J.J. Uhl}, 
``Vector Measures'', 
Math. Surveys, Vol. 15, Amer. Math. Soc., Providence, R.I.  (1977).

\smallskip\bibitem{Eng} R. Engelking, ``General Topology'', Revised edition, Heldermann Verlag, Berlin, 1989.

\smallskip\bibitem{EN} {\sc K.-J. Engel and R. Nagel}, 
``One-Parameter Semigroups for Linear Evolution Equations'',
Graduate Texts Math., Vol. 194, Springer-Verlag, 2000.

\smallskip\bibitem{fremlin} {\sc D. Fremlin, D. Garling and R. Haydon},
{\em Bounded measures on topological spaces},
Proc. London Math. Soc. {\bf 25} (1972), 115--136

\smallskip\bibitem{fu} {\sc M. Fuhrman}, 
{\em Analyticity of transition semigroups and closability of bilinear forms in
Hilbert spaces}, 
Studia Math. {\bf 115} (1995), 53--71.

\smallskip\bibitem{Go} {\sc J.A. Goldstein},
``Semigroups of Linear Operators and Applications'',  
Oxford University Press, Oxford, 1985.

\smallskip\bibitem{go} {\sc B. Goldys}, 
{\em On analyticity of Ornstein-Uhlenbeck semigroups}, 
Atti Accad. Naz. Lincei Cl. Sci. Fis. Mat. Natur. Rend. Lincei (9) Mat. Appl. 
{\bf 10} (1999), 131--140.

\smallskip\bibitem{clos} {\sc B. Goldys, F. Gozzi, and J.M.A.M. van Neerven}, 
{\em On closability of directional gradients}, Potential Anal. {\bf 18}  (2003), 289--310.

\smallskip\bibitem{GK} {\sc B. Goldys and M. Kocan}, 
{\em Diffusion semigroups in spaces of
continuous functions with mixed topology}, 
J. Differential Equations {\bf 173} (2001), 17--39.

\smallskip\bibitem{halmos} {\sc P. Halmos}, 
{\em Ten problems in Hilbert space},  
Bull. Amer. Math. Soc. {\bf 176}  (1970), 887--933.

\smallskip\bibitem{HPh} {\sc E. Hille and R.S. Phillips}, 
``Functional Analysis and Semi-Groups'',
Amer. Math. Soc. Colloq. Publ., Vol. XXXI, Rev. Ed., Providence, R.I., 1957
 

\smallskip\bibitem{lemerdy} {\sc C. Le Merdy},
{\em The similarity problem for bounded analytic semigroups on Hilbert spaces}, 
Semigroup Forum {\bf 56} (1998), 205--224.

\smallskip\bibitem{MR} {\sc Z.-M. Ma and M. R\"ockner}, 
``Introduction to the Theory of (Non-Symmetric) Dirichlet Forms'', 
Springer-Verlag, Berlin, 1992.

\smallskip\bibitem{Ne} {\sc J.M.A.M. van Neerven},
{\em Nonsymmetric Ornstein-Uhlenbeck semigroups in Banach spaces},
J. Funct. Anal. {\bf 155} (1998), 495--535.



\smallskip\bibitem{Ne3} {\sc  J.M.A.M. van Neerven}, 
{\em Uniqueness of invariant measures for the stochastic Cauchy problem in 
Banach spaces}, 
in: ``Recent Advances in Operator Theory and Related Topics: 
The Bela Sz\"okefalvi-Nagy Memorial Volume'', Operator theory: Advances and
Applications, Vol. 127, Birkhauser, 2001, pp. 491--517.

\smallskip\bibitem{NW} {\sc J.M.A.M. van Neerven and L. Weis}, 
{\em Weak limits and integrals of Gaussian covariances in Banach spaces},
Probab. Math. Statist. {\bf 25} (2005), 55--74 

\smallskip\bibitem{Pa} {\sc A. Pazy},
``Semigroups of Linear Operators and Applications to Partial Differential 
Equations'', Springer-Verlag, Berlin, 1983.

\smallskip\bibitem{priola}{\sc E. Priola},
{\em On a class of Markov type semigroups in spaces of uniformly 
continuous and bounded functions},  Studia Math.  {\bf 136}  (1999),  271--295.

\smallskip\bibitem{hiroshima} {\sc J. Pr\"uss and H. Sohr},
{\em Imaginary powers of elliptic second order differential operators in
$L^p$-spaces},
Hiroshima Math. J. {\bf 23} (1993), 161--192.


\smallskip\bibitem{Se} {\sc F.D. Sentilles}, {\em Bounded continuous 
functions on a completely regular space}, 
 Trans. Amer. Math. Soc.  {\bf 168}  (1972), 311--336.

\smallskip\bibitem{VTC} {\sc N.N. Vakhania, V.I. Tarieladze, and S.A. Chobanyan}, 
``Probability Distributions in Banach Spaces'',
D. Reidel Publishing Company, Dordrecht-Boston-Lancaster-Tokyo, 1987.

\smallskip\bibitem{Vu} {\sc V\~u Qu\^oc Ph\'ong}, 
{\em The operator equation $AX-XB = C$ with
unbounded operators $A$ and $B$ related to abstract Cauchy problems}, 
Math. Z. {\bf 208} (1991), 567--588.

\smallskip\bibitem{Wh} {\sc R.F. Wheeler}, 
{\em A survey of Baire measures and strict topologies}, 
Exposition. Math. {\bf 1} (1983),~97--190. 

\smallskip\bibitem{Wi} {\sc A.  Wiweger}, {\em Linear spaces with mixed topology}, 
Studia Math. {\bf 20} (1961), 47--68.

\smallskip\bibitem{Yo} {\sc S. Yosida}, ``Functional Analysis'', 
6th edition,  Grundlehren der Mathematischen Wissenschaften, Vol. 123, 
Springer-Verlag, Berlin-New York, 1980.

\smallskip\bibitem{Za} {\sc J. Zabczyk},
{\em Linear stochastic systems in Hilbert spaces; spectral
properties and limit behaviour}, 
Banach Center Publications {\bf 41}  (1985), 591--609.

\end{thebibliography}
\end{document}